\documentclass[11pt,a4paper]{amsart}
\usepackage{amsmath,amsthm,amssymb}
\usepackage[all]{xy}
\DeclareMathOperator{\sgn}{sgn}
\DeclareMathOperator{\tr}{tr}
\DeclareMathOperator{\pol}{pol}
\DeclareMathOperator{\opp}{opp}
\DeclareMathOperator{\Rea}{Re}
\DeclareMathOperator{\Pt}{Pt}
\DeclareMathOperator{\Sp}{Sp}
\DeclareMathOperator{\Step}{Step}
\DeclareMathOperator{\ad}{ad}
\DeclareMathOperator{\aad}{ad'}
\DeclareMathOperator{\Ad}{Ad}
\DeclareMathOperator{\cl}{cl}
\DeclareMathOperator{\Cl}{Cl}
\DeclareMathOperator{\Gr}{Gr}
\swapnumbers
\theoremstyle{definition}
\newtheorem{point}{}[section]
\newtheorem{remark}[point]{Remark}

\newtheorem{defi}[point]{Definition}
\newtheorem{example}[point]{Example}
\theoremstyle{plain}
\newtheorem{prop}[point]{Proposition}

\newtheorem{lemma}[point]{Lemma}
\newtheorem{theorem}[point]{Theorem}

\newtheorem{cor}[point]{Corollary}
\newcommand{\marginextend}[1]{ \addtolength{\oddsidemargin}{-#1}  \addtolength{\evensidemargin}{-#1}\addtolength{\textwidth}{#1}\addtolength{\textwidth}{#1}}
\newcommand{\updownextend}[1]{ \addtolength{\topmargin}{-#1}  \addtolength{\textheight}{#1}
\addtolength{\textheight}{#1}}
\marginextend{1.25cm}
\updownextend{0.0cm}
\begin{document}
\title{Fermionic quantum orthogonalizations I}
\date{\today}
\author{Gyula Lakos}
\email{lakos@cs.elte.hu}
\thanks{The author would like to thank the support of Bal\'azs Csik\'os.}
\address{Department of Geometry, E\"otv\"os University, P\'azm\'any P\'eter s.~1/C,  Budapest, H--1117, Hungary}
\begin{abstract}
We generalize classical orthogonalization procedures from real
linear algebra to the setting of fermionic quantum (FQ) operations.
In the case of the Gram-Schmidt orthogonalization procedure, the generalization is easy.
This, however, helps to obtain general information regarding FQ operations,
and to generalize the symmetric orthogonalization procedure.
\end{abstract}
\maketitle
\section*{Introduction}
In this paper, we study natural operations which assign to an $n$-tuple $(A_1,\ldots,A_n)$ of algebra-elements
another  $n$-tuple $(C_1,\ldots,C_n)$ of algebra-elements such that the identities
\[C_iC_j+C_jC_i=-2\delta_{ij}1\]
hold; i.~e., in our terminology, $(C_1,\ldots,C_n)$ forms a Clifford system.
(Due to Clifford (1879); see \cite{Lou} for reference and overview.)
Such operations generalize the classical orthogonalization procedures of linear algebra in the sense that if $V$ is a real vector space
equipped with a positive definite metric $q$, and $v_1,\ldots,v_n$ is an orthonormal set of vectors in $V$, then
the associated elements $\cl(v_1),\ldots,\cl(v_n)$ form a Clifford system in the associated Clifford algebra $\Cl(V,q)$;
which is a kind of fermionic quantized version of the original vector space.
(Fock (1928), Cook (1953); see \cite{BSZ} for references and about physical realizations;
except that our case would correspond to a space-time containing only finitely many points.)
For this reason we call such operations fermionic quantum (FQ) orthogonalization procedures.

Moreover, we also consider natural operations which assign to an $n$-tuple $(A_1,\ldots,A_n)$
another  $n$-tuple $(C_1,\ldots,C_n)$  such that  $(C_2C_1^{-1},\ldots,C_{n}C_1^{-1})$ is Clifford system;
i.~e., in our terminology, $(C_1,\ldots,C_n)$ forms a floating Clifford system.
(This condition is, in fact, not specific to index $1$ and to left or right.)
We call such operations fermionic quantum (FQ) conform-orthogonalization procedures, as they generalize
the usual conform-orthogonalization procedure. (Unfortunately,  mathematical terminology often uses
``orthogonalization'' instead of ``orthonormalization'', so we use ``conform-orthogonalization'' instead of ``orthogonalization''.)

This paper is organized as follows.
First, we will clarify what we mean by (natural analytic vectorial) FQ operations,
and, in particular, (conform-)orthogonalization procedures.
Then we will exhibit one such generalization of the Gram-Schmidt procedure.
This will help us to recover general information about formal FQ procedures, define a generalization  of the symmetric procedure
on the formal level, and then, extend the formal construction analytically.
\section{Fermionic quantum operations}
\begin{point}
In what follows, we will deal with unital real algebras $\mathfrak A$, such that

(i) $\mathfrak A$ is $\mathbb N$-filtered, $\mathfrak A=\mathfrak A^{(0)}\supset \mathfrak A^{(1)}\supset\ldots$, and  $\mathfrak A^{(i)}\mathfrak A^{(j)}\subset\mathfrak A^{(i+j)}$;

(ii) $\mathfrak A$ is complete with respect to the filtration; and Hausdorff, i.~e. $\bigcap_{i\in \mathbb N}\mathfrak A^{(i)}=0$;

(iii) all filtration quotients $\mathfrak A / \mathfrak A^{(i)}$ are finite-dimensional.

We assume $\mathfrak A$  is endowed with the natural topology induced from the filtration quotients
$\mathfrak A / \mathfrak A^{(i)}$, the $\infty$-topology.
In particular, a sequence converges in $\mathfrak A$, if it converges in every filtration quotient.
In certain cases it is also useful to consider the topology induced from the filtration quotient
$\mathfrak A / \mathfrak A^{(1)}$ only, the $0$-topology.
The set of invertible elements $\mathfrak A^\star$ in $\mathfrak A$ forms a $0$-open subset.
In what follows, let $\Theta$ denote the the natural factorization $\mathfrak A\rightarrow\mathfrak A / \mathfrak A^{(1)}$.

Such algebras are sufficiently interesting, but not that complicated analytically, so we call them ``convenient algebras''.
A principal example is the algebra $\mathfrak F_n$ generated by the symbols $Q_s$, $R_s$ $(s\in\{1,\ldots,n\})$,
and the relations $Q_iQ_j+Q_jQ_i=-2\delta_{ij}1$ with filtration induced from the grading $\deg Q_s=0$, $\deg R_s=1$.

Sometimes, in case of ambiguity, we use the notation $a\cdot b$ for the product of $a,b\in\mathfrak A$;
and we also use it in the situation when we take the product of an element of $\mathfrak A$ and an $n$-tuple from $\mathfrak A^n$ positionwise.
\end{point}

\begin{defi}
In what follows, an (analytic, natural, vectorial) FQ operation $\Psi$ will  mean an operation $\Psi$ satisfying the following  conditions
\texttt{(vR)}, \texttt{(Ana)}, \texttt{(Nat)}:

\texttt{(vR)} Vectorial domain and range condition: To an $n$-tuple $(A_1,\ldots,A_n)$ from
a convenient algebra $\mathfrak A$,
the operation  $\Psi$ may assign another $n$-tuple from $\mathfrak A$.

\texttt{(Nat)} Algebraic naturality:
(i) If  $\Psi(A_1,\ldots,A_n)=(B_1,\ldots,B_n)$ and $\phi:\mathfrak A\rightarrow \tilde{\mathfrak A}$ is a homomorphism, then
$\Psi(\phi(A_1),\ldots,\phi(A_n))=(\phi(B_1),\ldots,\phi(B_n))$.
In particular, if $\theta$ is invertible, then
$\Psi(\theta A_1\theta^{-1},\ldots,\theta A_n\theta^{-1})=(\theta B_1\theta^{-1},\ldots,\theta B_n\theta^{-1})$.
(ii) Moreover, we require that that domain of $\Psi$ on $\mathfrak A \oplus \mathfrak B$ should contain the
direct product of the domains on $\Psi$ on $\mathfrak A $  and $ \mathfrak B$.

\texttt{(Ana)} Analiticity:
If $Q=(Q_1,\ldots,Q_n)$ is a Clifford system, and $\Theta(R_1),\ldots, \Theta(R_n)\in \mathfrak A/\mathfrak A^{(1)} $ are sufficienty close to $0$, then
\[\Psi(Q_1+R_1,\ldots,Q_n+R_n)
=\left(\sum_{r=0}^\infty E^{\langle r\rangle}_{1},\ldots,\sum_{r=0}^\infty E^{\langle r\rangle}_{k},\ldots,\sum_{r=0}^\infty E^{\langle r\rangle}_{n}\right)
\tag{F1}\label{F1}\]
such that
\[E^{\langle r\rangle}_{k}=\sum_{\substack{\iota_0,\ldots,\iota_r\in \{0,1\}^n  ,\\ \, j_1,\ldots,j_r\in  \{1\ldots n\}}}
\psi_{\,\,\iota_0,\iota_1,\ldots,\iota_r}^{[k],j_1,\ldots,j_r}\,Q^{\iota_0}
R_{j_1}Q^{\iota_1}R_{j_2}\ldots R_{j_r}Q^{\iota_r},\tag{F2}\label{F2}\]
with some fixed data $\psi_{\,\,\iota_0,\iota_1,\ldots,\iota_r}^{[k],j_1,\ldots,j_r}\in\mathbb R$, independent from $Q$.
Here we  used the multipower notation, such that for $\iota=(\iota_1,\ldots,\iota_n)\in\{0,1\}^n$,
let $Q^{\iota}:=Q^{\iota_1}\cdot\ldots\cdot Q^{\iota_n}$.
(A  proper name would be ``local uniform analyticity''.
Nevetheless,  one can show that the  uniformity of the power series data $\psi_\star$
follows from \texttt{(Nat)}, by considering
appropriate elements in $\mathfrak A\oplus \mathfrak F_n $, anyway.)

Less fundamental but often considered conditions for FQ operations are:

\texttt{(vC)} Vectorial Clifford conservativity: If $(Q_1,\ldots,Q_n)$ is a Clifford system, then it yields $\Psi(Q_1,\ldots,Q_n)=(Q_1,\ldots,Q_n)$.

\texttt{(vL)}  Vectorial sign linearity: If  $\Psi(A_1,\ldots,A_n)=(B_1,\ldots,B_n)$ and $(\iota_1,\ldots,\iota_n)\in\{0,1\}^n$, then
$\Psi((-1)^{\iota_1}A_1,\ldots,(-1)^{\iota_n}A_n)=((-1)^{\iota_1}B_1,\ldots,(-1)^{\iota_n}B_n)$.

The FQ operation $\Psi$ is an FQ orthogonalization procedure if it also satisfies the additional condition

\texttt{(CP)} Clifford productivity: $\Psi$ yields Clifford systems.

A less fundamental but desirable condition for FQ orthogonalization procedures is

\texttt{($\mathtt{H^{0}}$)} Real scalar $0$-homogeneity: If  $\Psi(A_1,\ldots,A_n)=(B_1,\ldots,B_n)$ and $t>0$, then
$\Psi(tA_1,\ldots,tA_n)=(B_1,\ldots,B_n)$.
\end{defi}
\begin{remark}
One can also formulate (pseudo)scalar FQ operations with the (pseudo)scalar domain and range condition \texttt{(sR)$\equiv$(psR)} that $\Psi$ takes
values in $\mathfrak A$; and then \texttt{(Nat)} and \texttt{(Ana)} are understood accordingly.
The analogues of \texttt{(vC)}, \texttt{(vL)},  in the scalar case, are \texttt{(sC)} $\Psi(Q_1,\ldots,Q_n)=1$,  and
\texttt{(sL)} $\Psi((-1)^{\iota_1}A_1,\ldots,(-1)^{\iota_n}A_n)=\Psi(A_1,\ldots,A_n)$; and in the pseudoscalar case, the analogous conditions
are \texttt{(psC)} $\Psi(Q_1,\ldots,Q_n)=Q_1\ldots Q_n$, and  \texttt{(psL)}
$\Psi((-1)^{\iota_1}A_1,\ldots,(-1)^{\iota_n}A_n)=(-1)^{\iota_1+\ldots+\iota_n}\Psi(A_1,\ldots,A_n)$.
\end{remark}
The requirements for FQ conform-operations $\Psi$ are similar, but some conditions are modified.
\begin{defi}
In what follows, an (analytic natural vectorial) FQ conform-operation $\Psi$ will  mean an operation $\Psi$ satisfying the  conditions
\texttt{(vR)}, \texttt{(Nat)}, \texttt{(Biv')}, \texttt{(Ana')}; where $\texttt{(vR)}$ and \texttt{(Nat)} are as before, and the
other conditions are

\texttt{(Biv')} Bivariance: If  $\Psi(A_1,\ldots,A_n)=(B_1,\ldots,B_n)$ and  $\theta_1,\theta_2$ are invertible elements in $\mathfrak A$, then
$\Psi(\theta_1A_1\theta_2,\ldots,\theta_1A_n\theta_2)=(\theta_1B_1\theta_2,\ldots,\theta_1 B_n\theta_2)$.

\texttt{(Ana')} Floating analiticity:
If $Q=(Q_1,\ldots,Q_n)$ is a floating Clifford system, and $\Theta(R_1),\ldots, \Theta(R_n)\in \mathfrak A/\mathfrak A^{(1)} $ are sufficienty close to $0$, then we have an expansion as in \eqref{F1}, but with

\[E^{\langle r\rangle}_{k}=\sum_{\substack{\iota_0,\ldots,\iota_r\in \{0,1\}^{n-1}  ,\\ \, j_1,\ldots,j_r\in  \{1\ldots n\}}} \psi_{\,\,\iota_0,\iota_1,\ldots,\iota_r}^{[k],j_1,\ldots,j_r}\,
\tilde Q^{\iota_0}\tilde R_{j_1}\tilde Q^{\iota_1}\tilde R_{j_2}\ldots \tilde R_{j_r}\tilde Q^{\iota_r}Q_1,\tag{F2'}\label{F2'}\]
with some $\psi_{\,\,\iota_0,\iota_1,\ldots,\iota_r}^{[k],j_1,\ldots,j_r}\in\mathbb R$ given, independent from $Q$; where $\tilde R_k=R_kQ_1^{-1}$,
and $\tilde Q=(Q_2Q_1^{-1},\ldots,Q_nQ_1^{-1})$.
(The formula is  inspired by bivariance and taking a right translate by $Q_1^{-1}$, but it is not specific to index $1$
and to left or right.)

Less fundamental but generally considered conditions are
\texttt{(vC')} and \texttt{(vL)}; with

\texttt{(vC')} Vectorial floating Clifford conservativity: If $(Q_1,\ldots,Q_n)$ is a floating Clifford system, then $\Psi(Q_1,\ldots,Q_n)=(Q_1,\ldots,Q_n)$.

The FQ conform-operation $\Psi$ is called an FQ conform-orthogonalization procedure if it also satisfies the additional condition

\texttt{(CP')} Clifford productivity: $\Psi$ yields floating Clifford systems.
\end{defi}
\begin{remark}
One can define antivariant  FQ conform-operations by using

\texttt{(Biv'')} Antivariance: If  $\Psi(A_1,\ldots,A_n)=(B_1,\ldots,B_n)$ and
$\theta_1,\theta_2$ are invertible elements in $\mathfrak A$, then
$\Psi(\theta_1A_1\theta_2,\ldots,\theta_1A_n\theta_2)=(\theta_2^{-1}B_1\theta_1^{-1},\ldots, \theta_2^{-1}B_n\theta_1^{-1})$.

Then we need a modified form \texttt{(Ana'')} of analiticity: In  the right side of \eqref{F2'} $Q_1$ should be taken from the right
and $Q_1^{-1}$ should put the left.

In general, if $Q=(Q_1,\ldots,Q_n)$ is a floating Clifford system, then its inverse system $Q^{-1}=(Q_1^{-1},\ldots,Q_n^{-1})$
is also a floating Clifford system.
In particular, taking inverse systems of FQ conform-orthogonalizations yields antivariant FQ conform-operations.
\end{remark}

Even for real vector spaces there are various orthogonalization procedures.
Assume that $V$ is a real vector space equipped with a positive definite metric $q$.

\begin{point}
The most widely known orthogonalization method is the Gram-Schmidt procedure.
(Properly attributed to Gram (1879) and Schmidt (1907);
see \cite{LBG} for a review.)
First, the  Gram-Schmidt procedure gives the quasi-conform-orthogonal system
\[\tilde a_1=a_1, \quad \tilde a_2= \frac{\left| \begin{matrix}q(a_1,a_1)&a_1\\q(a_1,a_2)&a_2\end{matrix} \right|}
{q(a_1,a_1)},\quad \ldots,\quad\tilde a_n= \frac{\left|
\begin{matrix}q(a_1,a_1)&\cdots&q(a_{n-1},a_1)&a_1\\
\vdots&\ddots&\vdots&\vdots\\q(a_1,a_n)&\cdots&q(a_{n-1},a_n)&a_n\end{matrix} \right|}{\left|
\begin{matrix}q(a_1,a_1)&\cdots&q(a_{n-1},a_1)\\
\vdots&\ddots&\vdots\\q(a_1,a_{n-1})&\cdots&q(a_{n-1},a_{n-1})\end{matrix} \right|}.\]
Then we obtain the conform-orthogonal system
\[\check a_1=\sqrt{\frac{q(\tilde a_1,\tilde a_1)}{q(\tilde a_1,\tilde a_1)}}\,\tilde a_1,\quad\ldots,\quad\check
 a_n=\sqrt{\frac{q(\tilde a_1,\tilde a_1)}{q(\tilde a_n,\tilde a_n)}}\,\tilde a_n;\]
and the orthogonal system
\[\hat a_1=\frac{\tilde a_1}{\sqrt{q(\tilde a_1,\tilde a_1)}},\quad\ldots,\quad\hat a_n=\frac{\tilde a_n}{\sqrt{q(\tilde a_n,\tilde a_n)}}.\]

If one wants to generalize  the Gram-Schmidt (conform-)orthogonalization procedure, then it is reasonable require at least

\texttt{(PGS${}^($'${}^)$)} Permanence with respect to the Gram-Schmidt (conform-)orthogonalization procedure:
If  $\breve a_1,\ldots,\breve a_n$ is the Gram-Schmidt (conform-)orthogonalization of $a_1,\ldots,a_n$, then
$\Psi(\cl(a_1),\ldots,\cl(a_n))=(\cl(\breve a_1),\ldots,\cl(\breve a_n))$ in $\Cl(V,q)$.

\texttt{(Fil)} Filtration property: In the expression $\Psi(A_1,\ldots,A_n)=(B_1,\ldots,B_n)$, the term $B_i$ depends only on the terms $A_1,\ldots,A_i$
in the analytic expansion.

\texttt{(Par)} Parabolic property:
If $\Psi(A_1,\ldots,A_n)=(B_1,\ldots,B_n)$ and $U$ is a lower triangular matrix with positive elements in the diagonal, $U_{11}=1$,
then $\Psi(U(A_1,\ldots,A_n))=(B_1,\ldots,B_n)$.
\end{point}
\begin{point}
The symmetric orthogonalization procedure goes as follows.
(In pure mathematics, this procedure has no particular name, but it is
considered as a simple byproduct of the polar decomposition.
In this sense, it is to be associated to Autonne (1902) and Williamson (1935);
see the discussion in \cite{HJ}.
In quantum chemistry, it is known as  L\"owdin's symmetric orthogonalization procedure,
due to his work \cite{Low}, (1950); and, in some other areas, as the
Schweinler-Wigner method, after \cite{SW}, (1970);
both of these works refer back to the approximative solution by R. Landshoff (1936).)

Assume that the vectors $a_1,\dots,a_n$ form a non-degenerate system. In order to simplify the notation, let
\[S=\begin{bmatrix}q(a_1,a_1)&\cdots&q(a_{n},a_1)\\
\vdots&\ddots&\vdots\\q(a_1,a_{n})&\cdots&q(a_{n},a_{n})\end{bmatrix}.\]
Then the symmetric conform-orthogonalization procedure gives
$\check a_1,\ldots,\check a_n$, where
\[\begin{bmatrix}\check a_1&\ldots&\check a_n\end{bmatrix}=\begin{bmatrix}a_1&\ldots&a_n\end{bmatrix}
\sqrt{ S^{-1}}\cdot\frac{\tr\sqrt S}n;\]
and the symmetric orthogonalization procedure gives
$\hat a_1,\ldots,\hat a_n$, where
\[\begin{bmatrix}\hat a_1&\ldots&\hat a_n\end{bmatrix}=\begin{bmatrix}a_1&\ldots&a_n\end{bmatrix}
\sqrt{ S^{-1}}.\]
If one wants to generalize  the symmetric (conform-)orthogonalization procedure, then it is reasonable require at least

\texttt{(PSy${}^($'${}^)$)} Permanence with respect to the symmetric (conform-)orthogonalization procedure:
If  $\breve a_1,\ldots,\breve a_n$ is the symmetric (conform-)orthogonalization of $a_1,\ldots,a_n$, then
$\Psi(\cl(a_1),\ldots,\cl(a_n))=(\cl(\breve a_1),\ldots,\cl(\breve a_n))$ in $\Cl(V,q)$.

\texttt{($\mathtt\Sigma$)} Permutation symmetry: If $\Psi(A_1,\ldots,A_n)=(B_1,\ldots,B_n)$ and $\sigma\in\Sigma_n$ is a permutation,
then $\Psi(\sigma(A_1,\ldots,A_n))=\sigma(B_1,\ldots,B_n)$.

\texttt{(O)} Orthogonality: If $\Psi(A_1,\ldots,A_n)=(B_1,\ldots,B_n)$ and $U\in\mathrm O(n)$ is an orthogonal matrix,
then $\Psi(U(A_1,\ldots,A_n))=U(B_1,\ldots,B_n)$.

(Symmetry  is, of course,  weaker than orthogonality, but we prefer to include it here.)
\end{point}
\begin{point}
It is notable that both of the classical (conform-)orthogonalization procedures above satisfy,
at their level, the

\texttt{(FSt)} Fiber-star property: If $\Psi(A_1,\ldots,A_n)=(B_1,\ldots,B_n)$ and $t\in[0,1]$, then
$\Psi(tA_1+(1-t)B_1,\ldots,tA_n+(1-t)B_n)=(B_1,\ldots,B_n)$.
\end{point}

\section{On (floating) Clifford systems }
In our study, (floating) Clifford systems play the role of trivial objects.
Nevertheless, in order to deal with them, it is useful to introduce some terminology of geometric nature.
\begin{point}
If $R,Q\in\mathfrak A$, and $Q^2=-1$ or $Q^2=1$, then we set
\[R{}^{0}_{Q}:=\frac12 \left(R + Q^{-1}RQ\right),
\qquad%\]\[
R{}^{1}_{Q}:=\frac12 \left(R - Q^{-1}RQ\right);\]
these are  the symmetrization and antisymmetrization of $R$ with respect to $Q$.

We can apply this procedure multiple times, so we obtain  $R^{\iota_1}_{Q_1}\ldots{}^{\iota_s}_{Q_s}$.
It is easy to see that if $Q_1,\ldots,Q_s$ is a Clifford system, then the order of the (anti)symmetrizations can be interchanged.
Consequently, if $\sigma\in\Sigma_s$, then $R^{\iota_{\sigma(1)}}_{Q_{\sigma(1)}}\ldots{}^{\iota_{\sigma(s)}}_{Q_{\sigma(s)}}=R^{\iota_1}_{Q_1}\ldots{}^{\iota_s}_{Q_s}$.

\end{point}
\begin{point}
(a) If $R=(R_1,\ldots,R_n)$, $Q=(Q_1,\ldots,Q_n)$, and $Q$ is a Clifford system, $\iota=(\iota_1,\ldots,\iota_n)\in\{0,1\}^n $, then we define
\[\left(R / Q\right)^\iota_j:=\left(R_jQ_j^{-1}\right){}^{\iota_1}_{Q_1}\ldots^{\iota_n}_{Q_n},
\qquad%\]\[
\left(Q \backslash R\right)^\iota_j:=\left(Q_j^{-1}R_j\right){}^{\iota_1}_{Q_1}\ldots^{\iota_n}_{Q_n}.\]
These two terms are related according to the connecting property
\[ \left(Q \backslash R\right)^\iota_j
=Q_j^{-1} \left(R / Q\right)^\iota_j Q_j=(-1)^{\iota_j}
\left(R / Q\right)^\iota_j.\]

(b) Similarly, if $R=(R_1,\ldots,R_n)$, $Q=(Q_1,\ldots,Q_n)$, and $Q$ is a floating
Clifford system, $\iota=(\iota_1,\ldots,\iota_n)\in\{0,1\}^n $, then we define
\[\left(R / Q\right)^{\mathrm f \iota}_{j}:=\left(R_jQ_j^{-1}\right){}^{\iota_1-\iota_k}_{Q_1Q_k^{-1}}
\ldots^{\iota_n-\iota_k}_{Q_nQ_k^{-1}},
\qquad%\]\[
\left(Q \backslash R\right)^{\mathrm f\iota}_{j}:=\left(Q_j^{-1}R_j\right){}^{\iota_1-\iota_k}_{Q_k^{-1}Q_1}
\ldots^{\iota_n-\iota_k}_{Q_k^{-1}Q_n};\]
here $\iota_s-\iota_k$ is understood $\!\!\!\mod 2$, and
one can check that the definitions do not depend on the choice of $k$.
The two expressions above are related by the connecting property
\[Q_k^{-1}\cdot \left(R / Q\right)^{\mathrm f\iota}_{j}\cdot Q_k
=(-1)^{\iota_j-\iota_k} \left(Q \backslash R\right)^{\mathrm f\iota}_{ j}.\]
One can see that in the expressions above $\iota$ counts up to $\!\!\!\mod (1,\ldots,1)$; so the appropriate
index set is $\mathrm f\{0,1\}^n$, which is just $\{0,1\}^n$ but $\!\!\!\mod (1,\ldots,1)$.

(c) In the particular case when $Q$ is a Clifford system,
\[\left(R / Q\right)^{\iota}_{j}=\left(R / Q\right)^{\mathrm f \iota}_{j}+(-1)^{\iota_j}\left(Q\backslash R\right)^{\mathrm f \iota}_{j}.\]
\end{point}

\begin{point}
(a) Let $\mathfrak A^\star$ denote the set of units of $\mathfrak A$.
There are natural actions of $\mathfrak A^\star$ and $\mathfrak A$ on $\mathfrak A^n$ given by
\[(\Ad X)(A_1,\ldots,A_n)=(XA_1X^{-1},\ldots,XA_nX^{-1}),\]
and
\[(\ad X)(A_1,\ldots,A_n)=([X,A_1],\ldots,[X,A_n])\]
respectively. One has the exponential map $\exp:\mathfrak A\rightarrow \mathfrak A^\star$ as usual.

(b) Let $\mathfrak A^{\mathrm f}=\mathfrak A\times \mathfrak A^{\opp}$, with unit group $\mathfrak A^{\mathrm f\star}$.
Then $1^{\mathrm f}=(1,1^{\opp})\in\mathfrak A^{\mathrm f\star}$.
Then there are corresponding actions of $\mathfrak A^{\mathrm f\star}$ and $\mathfrak A^{\mathrm f}$ on $\mathfrak A^n$ given by
\[(\Ad^{\mathrm f}\, (X,Y^{\opp}))(A_1,\ldots,A_n)=(XA_1Y,\ldots,XA_nY),\]
and
\[(\ad^{\mathrm f} (X,Y^{\opp}))(A_1,\ldots,A_n)=(XA_1+A_1Y,\ldots,XA_n+A_nY),\]
respectively. Similarly, we have
$\exp^{\mathrm f}:\mathfrak A^{\mathrm f}\rightarrow \mathfrak A^{\mathrm f\star}$.

(c) One has the natural inclusions $\Delta:\mathfrak A\rightarrow\mathfrak A^{\mathrm f}$, $X\mapsto (X,-X^{\opp})$
and $\delta:\mathfrak A^\star\rightarrow\mathfrak A^{\mathrm f\star}$, $X\mapsto (X,(X^{-1})^{\opp})$
 compatible with the exponentials; etc.
\end{point}

\begin{point}
(a)
In what follows, let $\Gr\mathfrak A^n$ denote the subset Clifford systems of $\mathfrak A^n$.
If  $Q=(Q_1,\ldots,Q_n)$ is a Clifford system, then we define
\[\mathbf T_Q\Gr\mathfrak A^n:=\{(R_1,\ldots,R_n)\in\mathfrak A^n\,:\,\forall i,j\,\,  R_iQ_j+R_jQ_i+Q_iR_j+Q_jR_i=0\},\]
which can be rewritten as
\[\begin{split}&=\{(R_1,\ldots,R_n)\in\mathfrak A^n\,:\,\forall i\,\,(R_iQ_i^{-1})^0_{Q_i}=0;\,\,\forall i, j\,\, (R_iQ_i^{-1})^1_{Q_j}=(R_jQ_j^{-1})^1_{Q_i}\}\\
&=\{(R_1,\ldots,R_n)\in\mathfrak A^n\,:\,\exists X\,\, X{}^0_{Q_1}\ldots {}^0_{Q_n}=0 \,\,\&\,\,\forall i\,\,R_iQ_i^{-1}=2 X^1_{Q_i}\}\\
&=\{(R_1,\ldots,R_n)\in\mathfrak A^n\,:\,\exists X\,\, X{}^0_{Q_1}\ldots {}^0_{Q_n}=0 \,\,\&\,\,\forall i\,\,R_i=[X,Q_i]\}.\end{split}\]
Let
\[T_Q\mathfrak A:=\{X\in\mathfrak A\,:\, X{}^0_{Q_1}\ldots {}^0_{Q_n}=0\}.\]
It is easy to see the following:

(i) If $F:I\subset\mathbb R\rightarrow\Gr\mathfrak A^n$ is of class $C^1$ (on every filtration level), then $\dot F(t)\in\mathbf T_{F(t)}\Gr\mathfrak A^n$.

(ii) In particular, if $F:[0,1]\rightarrow\Gr\mathfrak A^n$ is given by
$F(t)=(\Ad \exp tX)Q$, then
$\dot F(0)=(\ad X)Q =([X,Q_1],\ldots,[X,Q_n])\in \mathbf T_{Q}\Gr\mathfrak A^n $.

(iii) The map
\[\aad Q :T_Q\Gr\mathfrak A^n\rightarrow\mathbf T_Q\Gr\mathfrak A^n\]
\[X\mapsto ([X,Q_1],\ldots,[X,Q_n])\]
gives a bijection.

For these reasons, it is safe to interpret $\mathbf T_Q\Gr\mathfrak A^n$ and $T_Q\Gr\mathfrak A^n$ as the tangent space
of $\Gr\mathfrak A^n$ at $Q$.

(b) Similarly, in what follows, let $\Gr^{\mathrm f}\mathfrak A^n$ denote the subset floating Clifford systems of $\mathfrak A^n$.
If we choose an index $k$, then  $Q=(Q_1,\ldots,Q_n)$ is a floating Clifford system if and only if
\[\forall i\neq k,j\neq k\,\,Q_iQ_k^{-1}Q_j+Q_jQ_k^{-1}Q_i=-2\delta_{i,j}Q_k\]
holds. We define
\begin{align*}
\mathbf T_Q\Gr^{\mathrm f}\mathfrak A^n:=&\{(R_1,\ldots,R_n)\in\mathfrak A^n\,:\,\forall i\neq k,j\neq k\,\,
R_iQ_k^{-1}Q_j-Q_iQ_k^{-1}R_kQ_k^{-1}Q_j+Q_iQ_k^{-1}R_j\\&
+R_jQ_k^{-1}Q_i-Q_jQ_k^{-1}R_kQ_k^{-1}Q_i+Q_jQ_k^{-1}R_i=-2\delta_{i,j}R_k\},
\end{align*}
which can be rewritten as
\[\begin{split}
&=\{(R_1,\ldots,R_n)\in\mathfrak A^n\,:\,\forall i\,\,(R_iQ_i^{-1}-(R_kQ_k^{-1})^0_{Q_iQ_k^{-1}})^0_{Q_iQ_k^{-1}}=0;\,\,
\\&\qquad\forall i, j\,\, (R_iQ_i^{-1}-(R_kQ_k^{-1})^0_{Q_iQ_k^{-1}})^1_{Q_jQ_k^{-1}}=(R_jQ_j^{-1}-(R_kQ_k^{-1})^0_{Q_jQ_k^{-1}})^1_{Q_iQ_k^{-1}}\}
\\&=\{(R_1,\ldots,R_n)\in\mathfrak A^n\,:\,\exists X,Y\,\, (X-Q_kYQ_k^{-1}){}^0_{Q_1Q_k^{-1}}\ldots {}^0_{Q_nQ_k^{-1}}=0 \,\,\&
\\&\qquad\,\,\forall i\,\,R_iQ_i^{-1}-(R_kQ_k^{-1})^0_{Q_iQ_k^{-1}}= (X-Q_kYQ_k^{-1})^1_{Q_iQ_k^{-1}} \,\,\&\,\, R_kQ_k^{-1}=(X+Q_kYQ_k^{-1}) \}
\\&=\{(R_1,\ldots,R_n)\in\mathfrak A^n\,:\,\exists X,Y\,\, (X-Q_kYQ_k^{-1}){}^0_{Q_1Q_k^{-1}}\ldots {}^0_{Q_nQ_k^{-1}}=0 \,\,\&\,\, R_i=XQ_i+Q_iY\}.\end{split}\]
Set
\[T_Q\Gr^{\mathrm f}\mathfrak A=\left\{(X,Y^{\opp})\in\mathfrak A\times\mathfrak A^{\opp}\,:\,(X-Q_kYQ_k^{-1}){}^{0}_{Q_1Q_k^{-1}}\ldots^{0}_{Q_nQ_k^{-1}}=0\right\}. \]
One can show that the definitions above do not depend on the choice of $k$.

(i') If $G:I\subset \mathbb R\rightarrow\Gr^{\mathrm f}\mathfrak A^n$ is of class $C^1$ (on every filtration level), then $\dot G(t)\in\mathbf T_{G(t)}\Gr^{\mathrm f}\mathfrak A^n$.

(ii') In particular, if  $G:[0,1]\rightarrow\Gr^{\mathrm f}\mathfrak A^n$ is given by $G(t)=
(\Ad^{\mathrm f} \exp^{\mathrm f}t(X,Y^{\opp}))Q$, then
$\dot G(0)= (\ad^{\mathrm f} (X,Y^{\opp}))Q=(XQ_1+Q_1Y,\ldots,XQ_n+Q_nY)\in \mathbf T_{Q}\Gr^{\mathrm f}\mathfrak A^n $.

(iii') The map
\[\aad^{\mathrm f} Q :T_Q\Gr^{\mathrm f}\mathfrak A^n\rightarrow\mathbf T_Q\Gr^{\mathrm f}\mathfrak A^n\]
\[(X,Y^{\opp})\mapsto (XQ_1+Q_1Y,\ldots,XQ_n+Q_nY)\]
gives a bijection.

For these reasons, it is safe to interpret $\mathbf T_Q\Gr^{\mathrm f}\mathfrak A^n$ and $T_Q\Gr^{\mathrm f}\mathfrak A^n$ as the tangent space
of $\Gr^{\mathrm f}\mathfrak A^n$ at $Q$.

(c)
Assume that $Q\in\Gr\mathfrak A^n$. Then the diagram
\[\xymatrix{ T_Q\Gr\mathfrak A^n\ar[r]^{\aad Q}\ar@/_/[d]_{\Delta_Q}
&\mathbf T_Q\Gr\mathfrak A^n\ar@/_/[d]_{\boldsymbol\Delta_Q}\\
T_Q\Gr^{\mathrm f}\mathfrak A^n\ar[r]_{\aad^{\mathrm f} Q}\ar@/_/[u]_{\nabla_Q}
&\mathbf T_Q\Gr^{\mathrm f}\mathfrak A^n\ar@/_/[u]_{\boldsymbol\nabla_Q}   }\]
with maps
\[\Delta_Q:\quad X\mapsto (X,-X^{\opp}),\qquad \boldsymbol\Delta_Q:\quad(X_1,\ldots,X_n)\mapsto(X_1,\ldots,X_n),\]
\[\nabla_Q:\quad(X,Y^{\opp})\mapsto\frac12(X-Y),\qquad \boldsymbol\nabla_Q:\quad (X_1,\ldots,X_n)\mapsto((X_1){}^1_{Q_1},\ldots,(X_n){}^1_{Q_n}),\]
shows how to split off the Clifford tangent space from the floating Clifford tangent space.
\end{point}
\begin{prop}
(a) If $F:[a,b]\rightarrow\Gr\mathfrak A^n$ is of class $C^1$ (in every filtration quotient),
then there is a unique $H:[a,b]\rightarrow\mathfrak A^\star$ of class $C^1$ such that
$H(a)=1$ and $\dot H(t)H(t)^{-1}=(\ad F(t))^{-1}\dot F(t)$ hold.
In this case, it yields   $ F(t)=(\Ad H(t) )F(a)$.

(b) If $G:[a,b]\rightarrow\Gr^{\mathrm f}\mathfrak A^n$ is of class $C^1$ (in every filtration quotient),
then there is a unique $K:[a,b]\rightarrow\mathfrak A^{\mathrm f\star}$ of class $C^1$ such that
$K(a)=1^{\mathrm f}$ and $ \dot K(t)K(t)^{-1}=(\ad^{\mathrm f} G(t))^{-1}\dot G(t)$ hold.
In this case, it yields $G(t)=(\Ad^{\mathrm f} K(t)) G(a)$.

(c) The latter construction extends the former one through the natural maps $\Delta, \delta$.
\begin{proof}
These are consequences of properties seen in the previous paragraph.
\end{proof}
\end{prop}
\begin{point}
Then it is safe to interpret the map $(\ad Q)^{-1} :\mathbf T_Q\Gr\mathfrak A^n\rightarrow T_Q\Gr\mathfrak A^n$
as the natural or minimal connection on the tangent space of $\Gr\mathfrak A^n$ at $Q$; and similarly in the floating case.
In the situation of the previous proposition, we set
\[\Pt(F):=H(b),\qquad\Pt^{\mathrm f}(G):=K(b),\]
respctively; as the parallel transport maps associated to $F$ and $G$.
\end{point}

\begin{point}\label{po:connec}
A splitting of the natural inclusion  $\mathbf T_Q\Gr\mathfrak A^n\hookrightarrow\mathfrak A^n$, or rather the corresponding map
$T_Q\Gr\mathfrak A^n\hookrightarrow\mathfrak A^n$, is called a connection over $Q\in\Gr\mathfrak A^n$.
Then, $\Pi_Q:\mathfrak A^n\rightarrow T_Q\Gr\mathfrak A^n$ being a connection is expressed by the identity that for
$X\in T_Q\Gr\mathfrak A^n$
\[\Pi_Q([X,Q_1],\ldots,[X,Q_n])=X.\tag{Conn}\label{C} \]

We will be interested in natural vectorial sign-linear connections which are as follows.
Let $\omega$ be a collection of real numbers  $\omega^{(\iota_1,\ldots,\iota_k)}_j$, $(j\in\{1,\ldots,n\},(\iota_1,\ldots,\iota_k)\in\{0,1\}^n)$ such that
\[\omega^{(\iota_1,\ldots,\iota_n)}_j=0\qquad\text{if}\quad \iota_j=0 \tag{CL1}\label{CL1}\]
and
\[\sum_{j=1}^n\omega^{(\iota_1,\ldots,\iota_n)}_j=1 \qquad\text{if}\quad (\iota_1,\ldots,\iota_n)\neq(0,\ldots,0)\tag{CL2}\label{CL2}\]
hold. Then $\omega$ defines a linear connection $\Pi^{\omega}_Q$ by the prescription
\[\Pi^{\omega}_QR:=\frac12\sum_{j\in\{1,\ldots,n\},(\iota_1,\ldots,\iota_n)\in\{0,1\}^n}\omega^{(\iota_1,\ldots,\iota_n)}_j
\cdot(R/Q)_j^{(\iota_1,\ldots,\iota_n)}.\]

(The point is that $\Pi^{\omega}$ looks like the expansion term of first order of a scalar
FQ operation with \texttt{(sL)} and satisfying the connection condition.)
Various special cases are:

(i) The Gram-Schmidt connection is defined using the data
\[\mathrm{GS}^{(\iota_1,\ldots,\iota_n)}_j:=
\begin{cases}1&\qquad\text{if}\quad j=\min\{h:\iota_h=1\}
\\0&\qquad\text{otherwise.}\quad\end{cases}\]
In this case
\[\Pi^{\mathrm{GS}}_QR=\frac12
\left(R_{Q_1}^1+R_{Q_1}^0{}_{Q_2}^1+\ldots+R_{Q_1}^0{}_{Q_2}^0\ldots{}_{Q_{n-1}}^0{}_{Q_n}^1\right).\]

(ii) The symmetric connection is defined using the data
\[\mathrm{Sy}^{(\iota_1,\ldots,\iota_k)}_j:=\frac{\iota_j}{\iota_1+\ldots+\iota_n}.\]
A useful formula (but not of particular use for us now) is
\[\Pi^{\mathrm{Sy}}_QR=
\frac12 \int_{h=0}^1\left(\sum_{j=1}^n  [R_j,Q_j^{-1}]\right){}^{[h}_{Q_1}{}^{[h}_{Q_2}\ldots{}^{[h}_{Q_n} \,\frac{\mathrm dh}h;\]
where we used the notation $X^{[ h}_Q=X^0_Q+hX^1_Q$.

(iii) If $w=(w_1,\ldots,w_n)$ is a positive weight vector, then the $w$-connection is defined from the data
\[(w)^{(\iota_1,\ldots,\iota_k)}_j:=\frac{w_j\iota_j}{w_1\iota_1+\ldots+w_n\iota_n}.\]

Then $\Pi^{\mathrm{Sy}}_Q$ corresponds to the case of equal weights, and $\Pi^{\mathrm{GS}}_Q$ corresponds to
the limit when $w_i/w_{i+1}\rightarrow+\infty$.
In particular, the weights $\mathrm{GS}(t)=(1,t,\ldots,t^{n-1})$, $t\in(0,1)$
connect the Gram-Schmidt and symmetric cases.
In fact, (\underline{iii})  one can easily extend this construction to the appropriate compactification $\widehat{P\mathbb R_+^n}$ of the space of (projectivized) positive weights. (We can just consider the closure of the connection data.)

More generally, we could consider natural connection expressions but without \texttt{(sL)}.
We do not go into details but one special case is  follows:

(iv) If $\eta=[\eta^{ij}]$ is a positive definite $n\times n$ matrix, then consider the map
\[M_{\eta,Q}:  T_Q\Gr\mathfrak A^n\rightarrow T_Q\Gr\mathfrak A^n \]
\[ X \mapsto \sum_{i,j}\eta^{ij}[[X,Q_j],Q_i] .\]
We can define
\[\Pi^{\eta}_QR:=\left(M_{\eta,Q}\right)^{-1}\left(\sum_{i,j}\eta^{ij}[R_j,Q_i]\right).\]
Indeed, if we apply an orthogonal change of coordinates in $(Q_1,\ldots,Q_n)$ and $(R_1,\ldots,R_n)$
such that $\eta$ becomes diagonal, then we see that our situation is essentially the same as in the weighted case.
An advantage of this latter formalism is that its invariance properties are more transparent.
Similarly, (\underline{iv}) one can extend this construction to the appropriate compactification $\widehat{P\mathfrak X_+^n(\mathbb R)}$
of the space of conformal positive definite forms. (Again, we can consider the closure of the connection data.)

In the cases above, the connections are natural in the sense that identity
\[\Pi_QQ=0\tag{CN1}\label{CN1}\]
holds, and they are (not at $Q$, but as collections) homomorphism-invariant, in particular, conjugation invariant:
\[(\Ad X)(\Pi_QR)=\Pi_{(\Ad X)Q}\,((\Ad X)R)\tag{CN2}\label{CN2}\]
(i.~e. they behave well as two-variable scalar FQ operation).

(b) A splitting of the natural inclusion  $\mathbf T_Q\Gr^{\mathrm f}\mathfrak A^n\hookrightarrow\mathfrak A^n$, or the corresponding map
$T_Q\Gr^{\mathrm f}\mathfrak A^n\hookrightarrow\mathfrak A^n$ is called a connection over $Q\in\Gr^{\mathrm f}\mathfrak A^n$.
Then, $\Pi_Q^{\mathrm f}:\mathfrak A^n\rightarrow T_Q\Gr^{\mathrm f}\mathfrak A^n$ being a connection is expressed by the identity that for
$(X,Y^{\opp})\in T_Q\Gr^{\mathrm f}\mathfrak A^n$
\[\Pi^{\mathrm f}_Q(XQ_1+Q_1Y,\ldots,XQ_n+Q_nY)=(X,Y^{\opp}).\tag{Conn'}\label{C'}\]

Natural vectorial sign-linear connections can be produced as earlier.
Let $\omega$ be a data as in \eqref{CL1}--\eqref{CL2}.
Set
\[\left(R / Q\right)^{\mathrm f\omega}:=
\sum_{\substack{ j\in\{1,\ldots,n\},(\iota_1,\ldots,\iota_n)\in \{0,1\}^n}}
\omega^{(\iota_1,\ldots,\iota_n)}_{j}\left(R / Q\right)^{\mathrm f(\iota_1,\ldots,\iota_n)}_{j};\]
and we define $\left(  Q\backslash R \right)^{\mathrm f\omega}$ analogously. Then let
\[\Pi^{\mathrm f\omega}_QR:=\left(\frac12\left(R / Q\right)^{\mathrm f\omega},\frac12\left(  Q\backslash R \right)^{\mathrm f\omega}\,{}^{\opp}\right).\]

Using exactly the same numerics as before we can define (i') $\Pi^{\mathrm f\mathrm{GS}}_Q$, (ii') $\Pi^{\mathrm f\mathrm{Sy}}_Q$,
(iii')$\Pi^{\mathrm fw}_Q$. Furthermore, we have analogously,

(iv') If $\eta=[\eta^{ij}]$ is a positive definite $n\times n$ matrix, then consider the map
\[M_{\eta,Q}^{\mathrm f}:  T_Q\Gr^{\mathrm f}\mathfrak A^n\rightarrow T_Q\Gr^{\mathrm f}\mathfrak A^n \]
\[ (X,Y^{\opp}) \mapsto \sum_{i,j}\eta^{ij} \left((XQ_j+Q_jY)Q_i^{-1},Q_i^{-1}(XQ_j+Q_jY)\,{}^{\opp}\right) .\]
We set
\[\Pi^{\mathrm f\eta}_QR:=\left(M_{\eta,Q}^{\mathrm f}\right)^{-1}\left(\sum_{i,j}\eta^{ij}\cdot(R_jQ_i^{-1},Q_i^{-1}R_j\,{}^{\opp})\right).\]

In the cases above, the connections are natural in the sense that identity
\[\Pi^{\mathrm f}_QQ=1^{\mathrm f}\tag{CN1'}\label{CN1'}\]
holds, and they are (not at $Q$, but as collections) homomorphism-invariant, and binvariant; in particular:
\[(\Ad^{\mathrm f} X)(\Pi^{\mathrm f} _QR)=\Pi^{\mathrm f} _{(\Ad^{\mathrm f} X)Q}\,((\Ad^{\mathrm f} X)R)
\tag{CN2'}\label{CN2'}\]
(i.~e. they behave well as two-variable scalar FQ conform-operations).

(c) If $Q$ is a Clifford system, then the connections above are related by
\[ \Pi_Q^{\omega}A=\nabla_Q\Pi_Q^{\mathrm f\omega}A.\]
\end{point}

\begin{lemma} \label{lem:conn}
If $Q=(Q_1,\ldots,Q_n)$ is a (floating) Clifford system, and $A=(A_1,\ldots,A_n)$, then

(i) $\Pi_Q^{\mathrm{GS}}A=0$ if and only if $(A_k)^1_{Q_1}{}^1_{Q_2}\ldots{}^1_{Q_{k-1}}{}^1_{Q_{k}}=0$ holds for all $k\in \{1,\ldots,n\}$.

(i') $\Pi_Q^{\mathrm{fGS}}A=1^{\mathrm f}$ if and only if $A_1=Q_1$ and $(A_k){}^1_{Q_2Q_1^{-1}}\ldots{}^1_{Q_{k-1}Q_1^{-1}}{}^1_{Q_{k}Q_1^{-1}}=0$ holds for all $k\in \{2,\ldots,n\}$.

(ii) $\Pi_Q^{\mathrm{Sy}}A=0$ if and only if $\sum_{i=1}^n[A_i,Q_i]=0$.

(ii') $\Pi_Q^{\mathrm{fSy}}A=1^{\mathrm f}$ if and only if $\sum_{i=1}^nA_iQ_i^{-1}=\sum_{i=1}^n Q_i^{-1}A_i=n\cdot 1$.

(iii) $\Pi_Q^{w}A=0$ if and only if $\sum_{i=1}^nw_i[A_i,Q_i]=0$.

(iii') $\Pi_Q^{\mathrm{f}w}A=1^{\mathrm f}$ if and only if $\sum_{i=1}^nw_iA_iQ_i^{-1}=\sum_{i=1}^n w_iQ_i^{-1}A_i=\left(\sum_{i=1}^n w_i\right)1$.

(iv) $\Pi_Q^{\eta}A=0$ if and only if $\sum_{i,j=1}^n\eta^{ij}[A_j,Q_i]=0$.

(iv') $\Pi_Q^{\mathrm{f}\eta}A=1^{\mathrm f}$ if and only if $\sum_{i,j=1}^n\eta^{ij}A_jQ_i^{-1}=\sum_{i,j=1}^n\eta^{ij}Q_i^{-1}A_j=\left(\sum_{i=1}^n\eta^{ii}\right)1$.
\begin{proof} (i) $(A_k)^1_{Q_1}{}^1_{Q_2}\ldots{}^1_{Q_{k-1}}{}^1_{Q_{k}}=0$ holds if and only if
$(A_kQ^{-1}_k)^0_{Q_1}{}^0_{Q_2}\ldots{}^0_{Q_{k-1}}{}^1_{Q_{k}}=0$; these latter terms vanish, however, if
and only if their sum vanishes. Case (ii) is similar; and the other cases are quite straightforward after passing
to decompositions with respect to $Q$.
\end{proof}
\end{lemma}
If we want extend the lemma to connections with compactified weights,
then we have to combine conditions of filtration ($\sim$(i)) and sum ($\sim$(ii--iv)) type.

\section{The generalization of the Gram-Schmidt procedure}
\begin{point}\label{po:pol}
 For $S\in\mathfrak A$, we define
\[S^{-1/2}:=\int_{t=0}^{2\pi} \frac{1}{\cos^2 t+S\sin^2 t}\,\frac{\mathrm dt}{2\pi},\]
assumed that $(\cos^2 t+S\sin^2 t)^{-1}$ exists for any $t$ (and then it is automatically continuous).
The domain condition is equivalent to $\Sp S\subset\mathbb C\setminus (-\infty,0]$
(the spectrum can be defined as usual).
It is well known that this operation produces an element commuting with $S$, and with the property $(S^{-1/2})^2= S^{-1}$.
In fact, $X=S^{-1/2}$ is characterized by the properties (i) $X^2=S^{-1}$, (ii) $\Sp X\subset \{z\in\mathbb C\,:\,\Rea z>0\}$.

 For $H\in\mathfrak A$, we can define  its polarization by
\[\pol H:=H (-H^2)^{-1/2}=\int_{t=0}^{2\pi} \frac{H}{\cos^2 t-H^2\sin^2 t}\,\frac{\mathrm dt}{2\pi},\]
assumed its defined.
It turns out that this definition is equivalent to the definition
\[\pol H:=\int_{t=0}^{2\pi} \frac{-\sin t+H\cos t}{\cos t+H\sin t}\,\frac{\mathrm dt}{2\pi},\]
defined if and only if $(\cos t+H\sin t)^{-1}$ exists for any $t$.
The domain condition is equivalent to $\Sp H\subset\mathbb C\setminus \mathbb R$.
This operation produces a skew-involution commuting with $H$.
In fact, $Q=\pol H$ can be characterized by the properties (i) $Q^{2}=-1$, (ii)   $H$ commutes with $Q$,
(iii) $\Sp HQ^{-1}\subset\{z\in\mathbb C\,:\,\Rea z>0\} $.

An advantage of these definitions among many others, is that
they are invariant with respect to homomorphisms $\phi:\mathfrak A\rightarrow\tilde{\mathfrak A}$.
(See\cite{Rin}, \cite{Gan}, and the historical remarks in \cite{DS}, and \cite{Hig}
on the origins of analytic functions on operators; a standard material by now.)
\end{point}

\begin{lemma} Some well-known properties of the polarization operation are:

(a) If $a,b,c,d\in\mathbb R$, $ad-bc\neq0$, then
$\pol\frac{aH+b}{cH+d}=\sgn(ad-bc)\pol H$.

(b) If $X$ (anti)commutes with $H$, then $X$ (anti)commutes with $\pol H$.

(c) If $t\in[0,1]$, then $\pol \left(tH+ (1-t)\pol H\right)=\pol H  $. \qed
\end{lemma}

\begin{point}
One can check that if $S=1+T$ is sufficiently close to $1$, then
\[(1+T)^{-1/2}=\sum_{r=0}^\infty\begin{pmatrix}-1/2\\r\end{pmatrix}T^{r}.\]
Consequently, if $H=Q+R$ is sufficiently close to a skew-involution $Q$, then
\begin{align}
\pol\, (Q+R) &=(Q+R)(1-QR-RQ)^{-1/2}\notag\\
&=\sum_{r=0}^\infty\biggl(\begin{pmatrix}-1/2\\r\end{pmatrix}Q(-QR-RQ)^{r}+\begin{pmatrix}-1/2\\r-1\end{pmatrix}R(-QR-RQ)^{r-1}\biggr).\label{e:exp}
\end{align}
\end{point}

\begin{defi}\label{def:GS}
(a) We define the Gram-Schmidt orthogonalization $\mathcal O^{\mathrm{GS}}$ of $(A_1,\ldots,A_n)$ as the $n$-tuple
$(Q_1,\ldots,Q_n)$, where
\[Q_1=\pol A_1,\quad Q_2=\pol \,\,(A_2)^1_{Q_1},\quad Q_3=\pol \,\,(A_3)^1_{Q_1}{}^1_{Q_2},\quad Q_4=\pol \,\,(A_3)^1_{Q_1}{}^1_{Q_2}{}^1_{Q_3},
\quad \ldots  ,\]
are defined recursively. ($\mathcal O^{\mathrm{GS}}$ is defined if every $Q_i$ is well-defined.)

(b) We define the Gram-Schmidt conform-orthogonalization $\mathcal O^{\mathrm{fGS}}$ of $(A_1,\ldots,A_n)$ as
$(A_1,\tilde Q_2A_1\ldots,\tilde Q_nA_1)$, where $(\tilde Q_2,\ldots,\tilde Q_n)$ is the
Gram-Schmidt orthogonalization of the $(n-1)$-tuple $(A_2A_1^{-1},\ldots,A_nA_1^{-1})$.
\end{defi}
\begin{theorem}
(a) $\mathcal O^{\mathrm{GS}}$ is an FQ orthogonalization procedure (\texttt{(vR)}, \texttt{(Nat)}, \texttt{(Ana)}, \texttt{(CP)})
with \texttt{(vC)}, \texttt{(vL)}, \texttt{(H\,$\mathtt{{}^{0}}$)}, which
generalizes the Gram-Schmidt orthogonalization procedure (\texttt{(PGS)}, \texttt{(Fil)}, \texttt{(Par)}),
and satisfies the additional condition
\texttt{(FSt)}.

(b) $\mathcal O^{\mathrm{fGS}}$ is an FQ conform-orthogonalization procedure
(\texttt{(vR)}, \texttt{(Nat)}, \texttt{(Biv')}, \texttt{(Ana')}, \texttt{(CP')})
with \texttt{(vC')}, \texttt{(vL)},
which generalizes the Gram-Schmidt conform-orthogonalization procedure (\texttt{(PGS')}, \texttt{(Fil)}, \texttt{(Par)}),
and satisfies the additional condition \texttt{(FSt)}.

(c) Furthermore, $\mathcal O^{\mathrm{GS}}$ and $\mathcal O^{\mathrm{fGS}}$ are defined on $0$-open subsets of $\mathfrak A^n$.
\begin{proof}
(a) All the basic properties are immediate from  the properties of the polarization operation,  except \texttt{(Ana)},
which holds due to \eqref{e:exp}.
\texttt{(PGS)} and \texttt{(Fil)} are also easy to see. If $A_k'=t_1A_1+\ldots+t_kA_k$, $t_k>0$,  then, due to $(A_s)_{Q_s}^1=0$, we have
$\pol\, (A_k')^1_{Q_1}{}^1_{Q_2}\ldots{}^1_{Q_{k-1}}=\pol\, t_k\cdot(A_k)^1_{Q_1}{}^1_{Q_2}\ldots{}^1_{Q_{k-1}}=Q_k$.
This proves \texttt{(Par)}. \texttt{(FSt)} follows from the corresponding property of polarization.
(b) This is similar to the previous one.
(c) The $0$-open property follows from the construction.
\end{proof}
\end{theorem}
\begin{theorem} \label{the:charGS}
$\mathcal O^{\mathrm{GS}}$ and $\mathcal O^{\mathrm{fGS}}$  can be characterized as follows:

(a) $\mathcal O^{\mathrm{GS}}(A_1,\ldots,A_n)=(Q_1,\ldots,Q_n)$ if and only if the following conditions hold:

\quad\texttt{(CP)} $(Q_1,\ldots,Q_n)$ is a Clifford system;

\quad\texttt{(LGS)} $(A_k)^1_{Q_1}{}^1_{Q_2}\ldots{}^1_{Q_{k-1}}{}^1_{Q_{k}}=0$ holds for all $k\in \{1,\ldots,n\}$; i.~e. $\Pi_Q^{\mathrm{GS}}A=0$;

\quad\texttt{(NSp)} $\Sp A_kQ_k^{-1}\subset \{z\in\mathbb C\,:\,\Rea z>0 \}$.

(b) $\mathcal O^{\mathrm{fGS}}(A_1,\ldots,A_n)=(Q_1,\ldots,Q_n)$ if and only if the following conditions hold:

\quad\texttt{(CP')} $(Q_1,\ldots,Q_n)$ is a floating Clifford system;

\quad\texttt{(LGS')} $A_1Q_1^{-1}=1$ and $(A_k){}^1_{Q_2Q_1^{-1}}\ldots{}^1_{Q_{k-1}Q_1^{-1}}{}^1_{Q_{k}Q_1^{-1}}=0$ holds for all $k\in \{2,\ldots,n\}$;\\
 i.~e. $\Pi_Q^{\mathrm{fGS}}A=0$;

\quad\texttt{(NSp)} $\Sp A_kQ_k^{-1}\subset \{z\in\mathbb C\,:\,\Rea z>0 \}$.
\begin{proof}
This  follows from the characterization of the polarization operation.
\end{proof}
\end{theorem}
We can term the condition \texttt{(LGS)} as the linear Gram-Schmidt condition, and \texttt{(NSp)} as the naive spectral condition.

\begin{remark} Not only $\mathcal O^{\mathrm{GS}}$, but its derivatives can be expressed by integral formulas.
Indeed,  make variations $A_i\mapsto A_i+\varepsilon_i$:
Then $Q_1=\pol A_1$ varies by
\[\int_{t=0}^{2\pi} \frac{1}{\cos^2 t-A_1^2\sin^2 t}(\varepsilon_1\cos^2 t+A_1\varepsilon_1A_1\sin^2t)
\frac{1}{\cos^2 t-A_1^2\sin^2 t}\,\frac{\mathrm dt}{2\pi};\]
yielding the expression for the derivative of the first term of $\mathcal O^{\mathrm{GS}}(A)$.
This implies that  $(A_2)^{1}_{Q_1}=\frac12(A_2+Q_1A_2Q_1)$ varies by
\begin{align*}
&\frac12\biggl(\varepsilon_2+Q_1\varepsilon_2Q_1
+\int_{t=0}^{2\pi} \frac{1}{\cos^2 t-A_1^2\sin^2 t}(\varepsilon_1\cos^2 t+A_1\varepsilon_1A_1\sin^2t)
\frac{1}{\cos^2 t-A_1^2\sin^2 t}\,\frac{\mathrm dt}{2\pi}A_2Q_1\\
&+Q_1A_2\int_{t=0}^{2\pi} \frac{1}{\cos^2 t-A_1^2\sin^2 t}(\varepsilon_1\cos^2 t+A_1\varepsilon_1A_1\sin^2t)
\frac{1}{\cos^2 t-A_1^2\sin^2 t}\,\frac{\mathrm dt}{2\pi}\biggr).
\end{align*}
After this, one can compute the variation of $\pol (A_2)^{1}_{Q_1}$, and continue in that manner,
yielding derivatives for further terms of $\mathcal O^{\mathrm{GS}}(A)$.
Ultimately, with $\varepsilon=(\varepsilon_1,\ldots,\varepsilon_n)$, we obtain the derivative
\[\partial\mathcal O^{\mathrm{GS}}(A;\varepsilon),\]
an $n$-tuple, where every term is linear but noncommutative in $\varepsilon$;
in fact, a closed integral formula using the additional terms $A_1,\ldots,A_n$ and  $Q_1,\ldots,Q_n$.
\end{remark}
\begin{lemma} Suppose that $\mathcal O^{\mathrm{GS}}(A)=Q$.

(i) If $\varepsilon\in\mathfrak A^n$, then  $\partial\mathcal O^{\mathrm{GS}}(A;\varepsilon)=\partial\mathcal O^{\mathrm{GS}}(A;(\aad Q)\Pi^{\mathrm{GS}}_Q\varepsilon)$.

(ii) If $X\in\mathfrak A$, then $\partial\mathcal O^{\mathrm{GS}}(A;(\ad X)A)=(\ad X)Q$.
\begin{proof}
(i) follows from the fact that $\Pi^{\mathrm{GS}}_Q\varepsilon=0$ implies $\partial\mathcal O^{\mathrm{GS}}(A;\varepsilon)=0$.
(ii) follows from the conjugation invariance of $\mathcal O^{\mathrm{GS}}$.
\end{proof}
\end{lemma}
One can make similar statements regarding the GS conform-orthogonalization $\mathcal O^{\mathrm{fGS}}$.

We remark that a consequence of the existence of the GS (conform-)orthogonalization procedure is the following
\begin{prop}\label{prop:conj}
(a) For $Q,R\in\Gr\mathfrak A^n$, let
\[\Pt^{\mathrm{GS}}(R,Q):=\Pt(t\in[0,1]\mapsto \mathcal O^{\mathrm{GS}}((1-t)Q+tR)).\]
Then $\Pt^{\mathrm{GS}}(R,Q)$ is defined in a $0$-open subset of $\Gr\mathfrak A^n \times \Gr\mathfrak A^n$ containing the diagonal,
natural and analytic (in appropriate sense), and satisfies the identity
\[(\Ad \Pt^{\mathrm{GS}}(R,Q))Q=R. \]
In particular, Clifford systems sufficiently close to each other are conjugates.

(b) For $Q,R\in\Gr^{\mathrm f}\mathfrak A^n$, let
\[\Pt^{\mathrm{fGS}}(R,Q):=\Pt(t\in[0,1]\mapsto \mathcal O^{\mathrm{fGS}}((1-t)Q+tR)).\]
Then $\Pt^{\mathrm{fGS}}(R,Q)$ is defined in a $0$-open subset of $\Gr^{\mathrm f}\mathfrak A^n \times \Gr^{\mathrm f}\mathfrak A^n$ containing the diagonal,
natural and analytic (in appropriate sense), and satisfies the identity
\[(\Ad^{\mathrm f}\Pt^{\mathrm{fGS}}(R,Q) )Q=R. \]
In particular, floating Clifford systems sufficiently close to each other are translations of each other.
\begin{proof}
This follows from the construction of the parallel transport.
\end{proof}
\end{prop}

We see that the generalization of  the Gram-Schmidt orthogonalization procedure was really straightforward.
But, when we want to generalize the symmetric orthogonalization procedure, the situation is not so simple.
An ideal solution would be an  explicit formula as in \ref{po:pol} and Definition \ref{def:GS}.
Such formulas, however, are just not so easy to find (for the author, at least).
Short of that, a characterization like in Theorem \ref{the:charGS} would be sufficient.
But such a characterization is not so straightforward either.
It is not apparent how to generalize the global spectral condition.
There are candidates, but, in general, dealing with global spectral conditions is not easy.
We circumvent this problem by giving a formal solution first, which we extend to an analytic solution.
Hence we will obtain global extensions of the ordinary symmetric orthogonalization procedure.

\section{Formal FQ operations}
Even the case of vector spaces shows that a natural orthogonalization procedure cannot be defined everywhere
(because degenerate systems cannot be  orthogonalized naturally), and there will be issues around domain and analiticity.
In order to avoid those complications, one can consider formal FQ operations.
\begin{point}
A formal FQ operation is just like an ordinary one, but analiticity \texttt{(Ana)} is replaced by
an explicit prescription to the domain \texttt{(Frm)}, and a weaker version of analiticity \texttt{(\underline{Ana})}:

\texttt{(Frm)} Formality: $(A_1,\ldots,A_n)$ is in the domain of $\Psi$ if and only if there is a
Clifford system $(Q_1,\ldots,Q_n)$, and $R_1,\ldots, R_n\in \mathfrak A^{(1)} $, such that $A_1=Q_1+R_1,\ldots,A_n=Q_n+R_n$.

\texttt{(\underline{Ana})} Formal or perturbative analiticity:
If $(Q_1,\ldots,Q_n)$ is a Clifford system, and $R_1,\ldots, R_n\in \mathfrak A^{(1)} $, then we have an analytic expansion as in \eqref{F1}--\eqref{F2}.

So, a (natural vectorial) formal FQ operation must satisfy \texttt{(vR)}, \texttt{(Nat)}, \texttt{(Frm)}, \texttt{(\underline{Ana})}.
(And, in good cases, \texttt{(vC)}, \texttt{(vL)}.)
One can see that $(A_1,\ldots,A_n)$ is in the formal domain if and only if $(\theta(A_1),\ldots,\theta(A_n) )$ is a Clifford system.
Indeed, this is obviously a necessary condition, while its converse follows from the existence of  $\mathcal O^{\mathrm{GS}}(A)$.
If one has an analytic FQ operation $\Psi$,
then its restriction $\underline{\Psi}$ to the formal domain will yield a formal FQ operation.
A formal FQ operation $\Psi$ is a formal FQ orthogonalization procedure, if \texttt{(CP)} is satisfied.

We remark, however, that in the  form above, the list of requirements for formal FQ operations is redundant:
\begin{lemma} Assuming \texttt{(vR)} and \texttt{(Frm)}, conditions
\texttt{(\underline{Ana})} and  \texttt{(Nat)} are equivalent.
\begin{proof}
\texttt{(\underline{Ana})}$\Rightarrow$\texttt{(Nat)}:
According to the conditions the power series expansion data $\psi_{\star}$ completely determines the
value of the operation, and the power series expansions and direct sums are compatible with homomorphisms.
\texttt{(Nat)}$\Rightarrow$\texttt{(\underline{Ana})}:
Consider the value of our operation applied to $(Q_1+R_1,\dots,Q_n+R_n)$ in the free algebra $\mathfrak F_n$.
This naturally yields a power series expansion.
Applying naturality, we conclude that it holds universally.
\end{proof}
\end{lemma}
Due to the universality of the algebra $\mathfrak F_n$,
the overall point is that the study of formal FQ operations can be considered as the algebraic study of appropriate
$n$-tuples from $\mathfrak F_n$ (i.~e. formal power series expansions),  a kind of state-field correspondence.
\end{point}
\begin{point}
Even if we consider only formal FQ operations, consequences, or infinitesimal versions
of other conditions, naturally, remain of interest.
According this, further conditions may or may not be to modified. For example,
\texttt{(H\,$\mathtt{{}^{0}}$)} is to be rephrased:

\texttt{(\underline{H}\,$\mathtt{{}^0}$)} Formal scalar $0$-homogeneity: If  $\Psi(A_1,\ldots,A_n)=(B_1,\ldots,B_n)$ and
$t$ is a formal commutative variable with filtration degree $1$, then
$\Psi((1+t)A_1,\ldots,(1+t)A_n)=(B_1,\ldots,B_n)$ in the appropriately filtered algebra $\mathfrak A[[t]]$.

Among the further conditions, \texttt{($\mathtt\Sigma$)}, \texttt{(O)}, \texttt{(Fil)} needs no modification.
But \texttt{(PSy)}, \texttt{(PGS)}, \texttt{(Par)} do, which are left to to the reader.
\texttt{(FSt)} looks like to be rephrased, but it is not necessary; in fact, due to the better convergence properties, if it holds, then
it holds in a stronger form:

\texttt{(\underline{FSt})} Formal fiber (affine) linearity: If $\Psi(A_1,\ldots,A_n)=(B_1,\ldots,B_n)$ and $t\in\mathbb R$, then
$\Psi(tA_1+(1-t)B_1,\ldots,tA_n+(1-t)B_n)=(B_1,\ldots,B_n)$.

A general strategy is to replace global invariance conditions with local ones, and then pass to the formal analogue.
But one has some choices here.
For example, group invariance conditions can be rephrased in terms of derivations; as a consequence,
for example, in the rephrased condition \texttt{(\underline{H}\,$\mathtt{{}^0}$)}, $t^2=0$ can be assumed.
Moreover, passing to derivations, or formal group actions, may be beneficial even in those cases,
when its utilization not necessarily required, like in the case of \texttt{(O)}.
Hence, one may and must be somewhat opportunistic about
conditions beyond the definition of formal FQ orthogonalization procedures.

When it comes to formal FQ conform-operations, in addition to
\texttt{(vR)}, \texttt{(Nat)}, \texttt{(Biv')}, conditions
 \texttt{(Frm')}, \texttt{(\underline{Ana}')}
must be phrased for floating Clifford systems; and the situation is, in general, similar.
\end{point}
\begin{theorem} Suppose that $\Psi$ is a formal FQ operation satisfying

\texttt{(\underline{Biv}')}
Formal bivariance: If  $\Psi(A_1,\ldots,A_n)=(B_1,\ldots,B_n)$ and  $\theta_1,\theta_2\in 1+\mathfrak A^{(1)}$, then
$\Psi(\theta_1A_1\theta_2,\ldots,\theta_1A_n\theta_2)=(\theta_1B_1\theta_2,\ldots,\theta_1 B_n\theta_2)$.

Then, we claim,

(a) $\Psi$ extends uniquely to a  formal FQ conform-operation. If $\Psi$ satisfies \texttt{(vL)}, then so does the extension.

(b) In particular: Suppose that $A=(A_1,\ldots,A_n)$, and $Q=(Q_1,\ldots,Q_n)$ is a Clifford system, $A-Q\in(\mathfrak A^{(1)})^{n}$, and
$H$ is such that $H^{2}=-1$ and $H$ anticommutes with all the $Q_i$. Then
\[\Psi(A\cdot \exp(tH) )=\Psi(A)\cdot \exp(tH) \]
for any $t\in\mathbb R$.
This, of course,  includes the special cases  $\exp(tH)=H,-1$.
\begin{proof}
(b)
If $t$ is a formal commutative variable of degree $1$, then due to \texttt{(\underline{Biv}')}, we know that
$\Psi(\exp(\frac t2H)\cdot A\cdot \exp(\frac t2H) )= \exp(\frac t2H)\cdot\Psi(A)\cdot\exp(\frac t2H) $ holds.
Then the same holds for real $t$,
due to the analyticity of the expression (on every filtration level) and that
we can use $Q=\exp(\frac t2H)\cdot Q\cdot \exp(\frac t2H)$ for the expansion base invariably.
Then we can conclude our statement from conjugation invariance.

(a) Assume that $A=(A_1,\ldots,A_n)$, and $Q=(Q_1,\ldots,Q_n)$ is a floating Clifford system such that
 $A-Q\in(\mathfrak A^{(1)})^{n}$.
Let $k\in\{1,\ldots,n\}$ be arbitrary. Then one can check that the elements
\[ \begin{bmatrix}Q_i&\\&(-1)^{\delta_{ik}}Q_i\end{bmatrix}\begin{bmatrix}&Q_k^{-1}\\Q_k^{-1}&\end{bmatrix}\]
$i\in\{1,\ldots,n\}$ form a Clifford system. We can define the extension $\Psi^{\mathrm e}(A_1,\ldots,A_n)$ by
\[ \begin{bmatrix}1\\0\end{bmatrix}^\top \cdot\Psi\left(\ldots,\begin{bmatrix}A_i&\\&(-1)^{\delta_{ik}}A_i\end{bmatrix}\begin{bmatrix}&A_k^{-1}\\A_k^{-1}&\end{bmatrix},\ldots\right)\cdot
\begin{bmatrix}&A_k\\A_k&\end{bmatrix}\begin{bmatrix}1\\0\end{bmatrix}. \]
This is clearly bivariant. The statement that this is indeed an extension of $\Psi$ follows from part (b). Unicity follows from bivariance.
\end{proof}
\end{theorem}
Similar conclusion can be drawn regarding antivariant conform-operations.

At this point, we have only the formal restrictions of the Gram-Schmidt procedures at hand.
A consequence of Theorem \ref{the:charGS} is
\begin{cor}\label{cor:charGS}
The formal restrictions
$\underline{\mathcal O}^{\mathrm{GS}}$ and $\underline{\mathcal O}^{\mathrm{fGS}}$  can be characterized as follows:

(a) $\underline{\mathcal O}^{\mathrm{GS}}(A)=Q$ if and only if (i) $Q$ is Clifford systems , (ii) $\Pi^{\mathrm{GS}}_QA=0$, (iii) $A-Q\in(\mathfrak A^{(1)})^n$.

(b) $\underline{\mathcal O}^{\mathrm{fGS}}(A)=Q$ if and only if (i) $Q$ is a floating Clifford system, (ii) $\Pi^{\mathrm{fGS}}_QA=1^{\mathrm f}$, (iii) $A-Q\in(\mathfrak A^{(1)})^n$.
\qed
\end{cor}
\begin{point}
In order to understand the power series expansions yielding
FQ operations, let us rewrite \eqref{F2} using the expansion
\[R_j=\sum_{\iota\in\{0,1\}^n} \left(R/Q\right)^\iota_j Q_j,\]
where every component either commutes or anticommutes with any $Q_l$. This yields
\[E^{\langle r\rangle}_{k}=\sum_{\substack{\iota_1,\ldots,\iota_r,\varkappa\in \{0,1\}^n  ,\\ \, j_1,\ldots,j_r\in  \{1\ldots n\} }} p_{\,\,\varkappa,\iota_1,\ldots,\iota_r}^{[k],j_1,\ldots,j_r}\,
(R/Q)^{\iota_1}_{j_1}\ldots (R/Q)^{\iota_r}_{j_r} Q^\varkappa Q_k\tag{F20}\label{F20}\]
where $ p_{\,\,\varkappa,\iota_1,\ldots,\iota_r}^{[k],j_1,\ldots,j_r}\in\mathbb R$.
We remark that these expressions use $2^{n(r+1)}n^{r+1}$ real coefficients on the $r$th expansion level, just like
the original \eqref{F2}.

Now, the implication of \texttt{(vL)} is that the power series can be reduced to the contributions coming from
$\varkappa=(0,\ldots,0)$; then
\[E^{\langle r\rangle}_{k}=\left(\sum_{\substack{\iota_1,\ldots,\iota_r\in \{0,1\}^n  ,\\ \, j_1,\ldots,j_r\in  \{1\ldots n\} }} p^{[k]}{}_{\iota_1,\ldots,\iota_r}^{j_1,\ldots,j_r}
(R/Q)^{\iota_1}_{j_1}\ldots (R/Q)^{\iota_r}_{j_r} \right)Q_k.\tag{F2L}\label{F2L}\]
This leaves $2^{nr}n^{r+1}$ real coefficients on the $r$th expansion level.

The meaning of $\texttt{(vC)}$ is simply that
\[E^{\langle0\rangle}_{k}=Q_k\tag{F2C}\label{F2C}.\]

However, there is a requirement hidden in  \texttt{(\underline{Ana})}:
This is the independence of the power series expansion from the base point.
In fact, this requirement can also be interpreted as conjugation invariance
with respect to elements from $1+\mathfrak A^{(1)}$, hidden in  \texttt{(Nat)}:
Indeed, as power series expansions
are invariant with respect to conjugation and possible base Clifford systems are conjugates
of each other, as Proposition \ref{prop:conj} shows; base point invariance is equivalent to
conjugation-invariance with respect to a fixed base point.
At first sight, dealing with this seems quite unaccessible, but this is not so:
\end{point}
\begin{theorem}
If we have an arbitrary collection  $p_{\star|\mathrm{GS}}$ of real numbers
\[p_{\,\,\varkappa,\iota_1,\ldots,\iota_r}^{[k],j_1,\ldots,j_r}\]
(in special case: without $\varkappa$) for such index sequences that
\[\forall s\qquad j_s\neq\min\{h: \iota_s(h)=1\};\]
then it can be extended uniquely to a collection $p_{\star}$ of real numbers with arbitrary indices, such that it yields an FQ operation
$\Psi$
(and, in special case: satisfying \texttt{(vL)}).

The recipe for $\Psi(A)$ is as follows: Let $Q=\mathcal O^{\mathrm{GS}}(A)$, $A=Q+R$, and apply the expansions
\eqref{F1}/\eqref{F20}  (in special case: \eqref{F1}/\eqref{F2L}) with respect to $p_{\star|\mathrm{GS}}$.

This leaves $((n-1)2^{n}+1)^r2^{n}n$  (in special case: $((n-1)2^{n}+1)^rn$ real coefficients)
on the $r$th expansion level to be chosen freely.
\begin{proof} First we note that $\Psi(A)$ is well-defined.
Indeed, if $Q=\mathcal O^{\mathrm{GS}}(A)$, then the terms $(R/Q)^\iota_j$ where $j=\min\{h: \iota(h)=1\}$ vanish
(``GS gauge'').
This means that in the expansion it does not matter that we know only $p_{\star|\mathrm{GS}}$.
On the other hand, these are the only  conditions for the terms $(R/Q)^\iota_j$; and this extension process should
reconstruct any $\Psi$ with the indicated properties.
\end{proof}
\end{theorem}

\begin{example}\label{exa:gs}
Let $\Psi$ be a formal FQ operation with \texttt{(vL)}, $n=2$.
In order to simply the notation, we will use the auxiliary terms
\[r_1=(R/Q)^{(00)}_1,\quad r_2=(R/Q)^{(01)}_1,\quad r_3=(R/Q)^{(10)}_1,\quad r_4=(R/Q)^{(11)}_1,\]
\[r_5=(R/Q)^{(00)}_2,\quad r_6=(R/Q)^{(01)}_2,\quad r_7=(R/Q)^{(10)}_2,\quad r_8=(R/Q)^{(11)}_2.\]
Then
\[\Psi(Q+R)_s=\left(p^{[s]}+\sum_{i=1}^8 p^{[s]}_ir_i+\sum_{i,j=1}^8 p^{[s]}_{ij}r_ir_j+O((R/Q)^{\geq3}) \right)Q_s,\qquad s\in\{1,2\}. \tag{F${}_2$}\label{U}\]
We can collect (some of) the coefficients into the scalar matrices $\mathbf P^{[s]}_0=[ p^{[s]}]$, the row matrices
$\mathbf P^{[s]}_1=[ p^{[s]}_i]_{i=1}^8$, and the square matrices $\mathbf P^{[s]}_2=[ p^{[s]}_{ij}]_{i=1}^8{}_{j=1}^8$.

If $\Psi=\underline{\mathcal O}^{\mathrm{GS}}$, then one can check that
$\mathbf P^{[1]}_0=[1]$, $ \mathbf P^{[2]}_0=[1]$, and
\[\mathbf P^{[1]}_1=\left[ \begin {array}{cccccccc} \boxed{0}&\boxed{0}&1&1&\boxed{0}&0&\boxed{0}&\boxed{0}\end {array} \right],\]
\[\mathbf P^{[2]}_1=\left[ \begin {array}{cccccccc} \boxed{0}&\boxed{0}&0&1&\boxed{0}&1&\boxed{0}&\boxed{0}\end {array} \right],\]
\[\mathbf P^{[1]}_2=\frac12\left[ \begin{array}{cccccccc}
\boxed{0}&\boxed{0}&-1&-1&\boxed{0}&0&\boxed{0}&\boxed{0}\\
\boxed{0}&\boxed{0}&-1&-1&\boxed{0}&0&\boxed{0}&\boxed{0}\\
-1&-1&1&1&0&0&0&0\\
-1&-1&1&1&0&0&0&0\\
\boxed{0}&\boxed{0}&0&0&\boxed{0}&0&\boxed{0}&\boxed{0}\\
0&0&0&0&0&0&0&0\\
\boxed{0}&\boxed{0}&0&0&\boxed{0}&0&\boxed{0}&\boxed{0}\\
\boxed{0}&\boxed{0}&0&0&\boxed{0}&0&\boxed{0}&\boxed{0}
\end{array}\right],\]
\[\mathbf P^{[2]}_2=\frac12\left[\begin {array}{cccccccc}
\boxed{0}&\boxed{0}&0&-1&\boxed{0}&0&\boxed{0}&\boxed{0}\\
\boxed{0}&\boxed{0}&-1&0&\boxed{0}&0&\boxed{0}&\boxed{0}\\
0&-1&0&1&0&1&0&-1\\
-1&0&-1&1&0&1&-1&0\\
\boxed{0}&\boxed{0}&0&0&\boxed{0}&-1&\boxed{0}&\boxed{0}\\
0&0&-1&1&-1&1&0&0\\
\boxed{0}&\boxed{0}&0&-1&\boxed{0}&0&\boxed{0}&\boxed{0}\\
\boxed{0}&\boxed{0}&1&0&\boxed{0}&0&\boxed{0}&\boxed{0}
\end {array} \right].\]
Notice that all the boxed entries vanish.
Indeed, if we set $r_3\equiv(R/Q)^{(10)}_1=0$, $r_4=(R/Q)^{(11)}_1=0$, $r_6=(R/Q)^{(01)}_2=0$, i.~e. $\Pi_Q^{\mathrm{GS}}(Q+R)=0$,
then $\underline{\mathcal O}^{\mathrm{GS}}(Q+R)=Q$, so the pure $\{1,2,5,7,8\}$-terms must vanish in expansion degrees $r\geq1$.

In case $\Psi$ is arbitrary, we can compute   \eqref{U} using the base $Q$, i.~e.~as it is.
But alternatively, we can compute using  $\tilde Q=\mathcal O^{\mathrm{GS}}(A)$ as the base Clifford system.
Then, considering the result, we obtain
$\mathbf P^{[1]}_0=[\boxed{p^{[1]}}]$, $ \mathbf P^{[2]}_0=[\boxed{p^{[2]}}]$, and
\[\mathbf P^{[1]}_1=\left[ \begin {array}{cccccccc} \boxed{p_{{1}}^{[1]}}&\boxed{p_{{2}}^{[1]}}&p_
{{}}^{[1]}&p^{[1]}-p_{{8}}^{[1]}&\boxed{p_{{5}}^{[1]}}&0&\boxed{p_{{7}}^{[1]}}&
\boxed{p_{{8}}^{[1]}}\end {array} \right],\]
\[\mathbf P^{[2]}_1=\left[ \begin {array}{cccccccc} \boxed{p_{{1}}^{[2]}}&\boxed{p_{{2}}^{[2]}}&0&p^{[2]}-
p_{{8}}^{[2]}&\boxed{p_{{5}}^{[2]}}&p^{[2]}&\boxed{p_{{7}}^{[2]}}
&\boxed{p_{{8}}^{[2]}}\end {array} \right],\]
\[\mathbf P^{[1]}_2=\left[ \begin {smallmatrix}
\boxed{p_{{11}}^{[1]}}&\boxed{p_{{12}}^{[1]}}&-\frac12p^{[1]}+\frac12p_{{1}}^{[1]}\\
\boxed{p_{{21}}^{[1]}}&\boxed{p_{{22}}^{[1]}}&-\frac12p^{[1]}+\frac12p_{{2}}^{[1]}+\frac12p_{{8}}^{[1]}\\
-\frac12p^{[1]}+\frac12p_{{1}}^{[1]}&-\frac12p^{[1]}+\frac12p_{{2}}^{[1]}+\frac12p_{{8}}^{[1]}&\frac12p^{[1]}-\frac12p_{{1}}^{[1]}\\
-\frac12p^{[1]}+\frac12p_{{1}}^{[1]}+\frac12p_{{8}}^{[1]}-p_{{81}}^{[1]}&-\frac12p^{[1]}+\frac12p_{{2}}^{[1]}-p_{{82}}^{[1]}
&\frac12p^{[1]}-\frac12p_{{2}}^{[1]}-\frac12p_{{8}}^{[1]}\\
\boxed{p_{{51}}^{[1]}}&\boxed{p_{{52}}^{[1]}}&\frac12p_{{5}}^{[1]}+\frac12p_{{7}}^{[1]}\\
\frac12p_{{1}}^{[1]}-\frac12p_{{2}}^{[1]}&-\frac12p_{{1}}^{[1]}+\frac12p_{{2}}^{[1]}&\frac12p_{{8}}^{[1]}\\
\boxed{p_{{71}}^{[1]}}&\boxed{p_{{72}}^{[1]}}&\frac12p_{{5}}^{[1]}+\frac12p_{{7}}^{[1]}\\ \boxed{p_{{81}}^{[1]}}&\boxed{p_{{82}}^{[1]}}&\frac12p_{{8}}^{[1]}\end {smallmatrix} \right|\ldots\]
\[\ldots\left| \begin {smallmatrix}
-\frac12p^{[1]}+\frac12p_{{1}}^{[1]}+\frac12p_{{8}}^{[1]}-p_{{18}}^{[1]}&\boxed{p_{{15}}^{[1]}}&-\frac12
p_{{1}}^{[1]}+\frac12p_{{2}}^{[1]}&\boxed{p_{{17}}^{[1]}}&\boxed{p_{{18}}^{[1]}}\\
-\frac12p^{[1]}+\frac12p_{{2}}^{[1]}-p_{{28}}^{[1]}&\boxed{p_{{25}}^{[1]}}&\frac12p_{{1}}^{[1]}-\frac12p_{{2}}^{[1]}
&\boxed{p_{{27}}^{[1]}}&\boxed{p_{{28}}^{[1]}}\\
\frac12p^{[1]}-\frac12p_{{2}}^{[1]}-\frac12p_{{8}}^{[1]}&\frac12p_{{5}}^{[1]}-\frac12p_{{7}}^{[1]}&
-\frac12p_{{8}}^{[1]}&-\frac12p_{{5}}^{[1]}+\frac12p_{{7}}^{[1]}&\frac12p_{{8}}^{[1]}\\
\frac12p^{[1]}-\frac12p_{{1}}^{[1]}+\frac12p_{{5}}^{[1]}-p_{{8}}^{[1]}+p_{{88}}^{[1]}&\frac12p_{{5}}^{[1]}-\frac12p_{{8}}^{[1]}-p_{{85}}^{[1]}
&\frac12p_{{8}}^{[1]}&\frac12p_{{7}}^{[1]}-p_{{87}}^{[1]}&-\frac12p_{{5}}^{[1]}+\frac12p_{{8}}^{[1]}-p_{{88}}^{[1]}\\
\frac12p_{{5}}^{[1]}-\frac12p_{{8}}^{[1]}-p_{{58}}^{[1]}&\boxed{p_{{55}}^{[1]}}&-\frac12p_{{5}}^{[1]}&\boxed{p_{{57}}^{[1]}}&\boxed{p_{{58}}^{[1]}}\\
-\frac12p_{{8}}^{[1]}&\frac12p_{{5}}^{[1]}&-\frac12p_{{5}}^{[1]}&-\frac12p_{{8}}^{[1]}+\frac12p_{{7}}^{[1]}&-\frac12p_{{7}}^{[1]}+\frac12p_{{8}}^{[1]}\\
\frac12p_{{7}}^{[1]}-p_{{78}}^{[1]}&\boxed{p_{{75}}^{[1]}}&-\frac12p_{{7}}^{[1]}-\frac12p_{{8}}^{[1]}&\boxed{p_{{77}}^{[1]}}&\boxed{p_{{78}}^{[1]}}\\
-\frac12p_{{5}}^{[1]}+\frac12p_{{8}}^{[1]}-p_{{88}}^{[1]}&\boxed{p_{{85}}^{[1]}}
&-\frac12p_{{7}}^{[1]}-\frac12p_{{8}}^{[1]}&\boxed{p_{{87}}^{[1]}}&\boxed{p_{{88}}^{[1]}}\end {smallmatrix} \right],\]
\[\mathbf P^{[2]}_2=\left[ \begin {smallmatrix}
\boxed{p_{{11}}^{[2]}}&\boxed{p_{{12}}^{[2]}}&-\frac12p_{{1}}^{[2]}\\
\boxed{p_{{21}}^{[2]}}&\boxed{p_{{22}}^{[2]}}&-\frac12p^{[2]}-\frac12p_{{2}}^{[2]}+\frac12p_{{8}}^{[2]}\\
\frac12p_{{1}}^{[2]}&-\frac12p^{[2]}+\frac12p_{{2}}^{[2]}+\frac12p_{{8}}^{[2]}&-\frac12p_{{1}}^{[2]}\\
-\frac12p^{[2]}+\frac12p_{{1}}^{[2]}+\frac12p_{{8}}^{[2]}-p_{{81}}^{[2]}&\frac12p_{{2}}^{[2]}-p_{{82}}^{[2]}
&-\frac12p^{[2]}-\frac12p_{{2}}^{[2]}+\frac12p_{{8}}^{[2]}\\
\boxed{p_{{51}}^{[2]}}&\boxed{p_{{52}}^{[2]}}&-\frac12p_{{5}}^{[2]}+\frac12p_{{7}}^{[2]}\\
\frac12p_{{1}}^{[2]}-\frac12p_{{2}}^{[2]}&-\frac12p_{{1}}^{[2]}+\frac12p_{{2}}^{[2]}&-\frac12p^{[2]}+\frac12p_{{8}}^{[2]}\\
\boxed{p_{{71}}^{[2]}}&\boxed{p_{{72}}^{[2]}}&\frac12p_{{5}}^{[2]}-\frac12p_{{7}}^{[2]}\\
\boxed{p_{{81}}^{[2]}}&\boxed{p_{{82}}^{[2]}}&\frac12p^{[2]}-\frac12p_{{8}}^{[2]}
\end {smallmatrix} \right|\ldots
\]
\[\ldots\left| \begin {smallmatrix}
-\frac12p^{[2]}+\frac12p_{{1}}^{[2]}+\frac12p_{{8}}^{[2]}-p_{{18}}^{[2]}&\boxed{p_{{15}}^{[2]}}&\frac12p_{{1}}^{[2]}+\frac12
p_{{2}}^{[2]}&\boxed{p_{{17}}^{[2]}}&\boxed{p_{{18}}^{[2]}}\\
\frac12p_{{2}}^{[2]}-p_{{28}}^{[2]}&\boxed{p_{{25}}^{[2]}}&\frac12p_{{1}}^{[2]}+\frac12p_{{2}}^{[2]}&\boxed{p_{{27}}^{[2]}}&\boxed{p_{{28}}^{[2]}}\\
\frac12p^{[2]}-\frac12p_{{2}}^{[2]}-\frac12p_{{8}}^{[2]}&\frac12p_{{5}}^{[2]}-\frac12p_{{7}}^{[2]}&\frac12p^{[2]}-\frac12p_{{8}}^{[2]}&
-\frac12p_{{5}}^{[2]}+\frac12p_{{7}}^{[2]}&-\frac12p^{[2]}+\frac12p_{{8}}^{[2]}\\
\frac12p^{[2]}-\frac12p_{{1}}^{[2]}+\frac12p_{{5}}^{[2]}-p_{{8}}^{[2]}+p_{{88}}^{[2]}&\frac12p_{{5}}^{[2]}-\frac12p_{{8}}^{[2]}-p_{{85}}^{[2]}&
\frac12p^{[2]}-\frac12p_{{8}}^{[2]}&-\frac12p^{[2]}+\frac12p_{{7}}^{[2]}-p_{{87}}^{[2]}&-\frac12p_{{5}}^{[2]}+\frac12p_{{8}}^{[2]}-p_{{88}}^{[2]}\\
-\frac12p_{{8}}^{[2]}+\frac12p_{{5}}^{[2]}-p_{{58}}^{[2]}&\boxed{p_{{55}}^{[2]}}&-\frac12p^{[2]}+\frac12p_{{5}}^{[2]}
&\boxed{p_{{57}}^{[2]}}&\boxed{p_{{58}}^{[2]}}\\
\frac12p^{[2]}-\frac12p_{{8}}^{[2]}&-\frac12p^{[2]}+\frac12p_{{5}}^{[2]}&\frac12p^{[2]}-\frac12p_{{5}}^{[2]}&
\frac12p_{{7}}^{[2]}-\frac12p_{{8}}^{[2]}&-\frac12p_{{7}}^{[2]}+\frac12p_{{8}}^{[2]}\\
-\frac12p^{[2]}+\frac12p_{{7}}^{[2]}-p_{{78}}^{[2]}&\boxed{p_{{75}}^{[2]}}&\frac12p_{{7}}^{[2]}-\frac12p_{{8}}^{[2]}
&\boxed{p_{{77}}^{[2]}}&\boxed{p_{{78}}^{[2]}}\\
-\frac12p_{{5}}^{[2]}+\frac12p_{{8}}^{[2]}-p_{{88}}^{[2]}&\boxed{p_{{85}}^{[2]}}&
-\frac12p_{{7}}^{[2]}+\frac12p_{{8}}^{[2]}&\boxed{p_{{87}}^{[2]}}&\boxed{p_{{88}}^{[2]}}
\end {smallmatrix} \right].\]

What we see that we can eliminate the coefficients which are not in the boxed positions.
Indeed, expansion in the Gram-Schmidt base uses coefficients with pure $\{1,2,5,7,8\}$-indices.
We, however, cannot possibly get any further information about pure $\{1,2,5,7,8\}$-terms, because ``power series
expansion in Gram-Schmidt base'' is a natural formal FQ operation, and it yield different results if we set
$r_3=r_4=r_6=0$.

Reducing to pure $\{1,2,5,7,8\}$-indices makes the reduction not so symmetric, but
the collection of all the eliminating equations as a whole is certainly symmetric:
The conjugation-invariance property itself is symmetric for permutations in the variables $s\in\{1,2\}$, hence
there is a symmetry induced by
\[r_1=(R/Q)^{(00)}_1\leftrightarrow r_5=(R/Q)^{(00)}_2,\qquad r_2=(R/Q)^{(01)}_1\leftrightarrow r_7=(R/Q)^{(10)}_2,\]
\[r_3=(R/Q)^{(10)}_1\leftrightarrow r_6=(R/Q)^{(01)}_2,\qquad r_4=(R/Q)^{(11)}_1 \leftrightarrow r_8=(R/Q)^{(11)}_2.\]
Considering the permutation \[\varpi=(15)(27)(36)(48),\] we see that the set of equations is symmetric
with respect to $p^{[1]}_{i_1, \ldots } \leftrightarrow p^{[2]}_{\varpi(i_1), \ldots } $.

For example,
\[p_{{36}}^{[1]}= -\frac12p_{{8}}^{[1]}\]
yields
\[p_{{63}}^{[2]}= -\frac12p_{{4}}^{[2]};\]
which is, we know, valid, due to
\[p_{{63}}^{[2]}=-\frac12p^{[2]}+\frac12 p_{{8}}^{[2]}  \quad\text{and}\quad  p_{{4}}^{[2]}=p^{[2]}-p_{{8}}^{[2]}.\]
\end{example}
In particular, there are plenty of FQ operations. The situation is similar with respect to FQ orthogonalization procedures:
\begin{theorem}
If we have an arbitrary collection $\tilde p_{\star|\mathrm{GS}}$ of real numbers
\[\tilde p_{\iota_1,\ldots,\iota_r}^{j_1,\ldots,j_r}\]
for such index sequences that
\[r\geq 1,\qquad \forall s\quad j_s\neq\min\{h: \iota_s(h)=1\}, \qquad\sum_{s=1}^r\iota_s\neq(0,\ldots,0);\]
then it  yields a formal FQ orthogonalization procedure
$\Psi$ satisfying \texttt{(vL)}, \texttt{(vC)}  with the prescription
\[(\Ad \exp E^{\langle 1\rangle})(\Ad\exp E^{\langle 2\rangle})(\Ad\exp E^{\langle 3\rangle})\ldots\mathcal O^{\mathrm{GS}}(A).\]
where $Q=\mathcal O^{\mathrm{GS}}(A)$, $A=Q+R$, and
\[E^{\langle r\rangle}=\sum_{\substack{\iota_1,\ldots,\iota_r\in \{0,1\}^n  ,\\ \, j_1,\ldots,j_r\in  \{1\ldots n\} }}
\tilde p_{\iota_1,\ldots,\iota_r}^{j_1,\ldots,j_r}
(R/Q)^{\iota_1}_{j_1}\ldots (R/Q)^{\iota_r}_{j_r} .\]
This assignment gives a bijection between such data $\tilde p_{\star|\mathrm{GS}}$ and such FQ orthogonalization procedures.

This leaves $(((n-1)2^{n}+1)^r-1)(1-2^{-n})$ real coefficients on the $r$th expansion level $(r\geq1)$ to be chosen freely.
\begin{proof}
According to Proposition \ref{prop:conj}, such formal FQ orthogonalization procedures are conjugate to $\underline{\mathcal O}^{\mathrm{GS}}$.
Writing the conjugating terms in an appropriate form, we can achieve a most economical form as above.
\end{proof}
\end{theorem}
\begin{point}
We can proceed in the conform case similarly.
For example, if $R=(R_1,\ldots,R_n)$, $Q=(Q_1,\ldots,Q_n)$, and $Q$ forms a floating Clifford system,
then one can show that  \texttt{(vL)} implies that the operation $\Psi$ allows an expansion around $Q$ as in \eqref{F1}, but with
\[E^{\langle r\rangle}_{k}=\left(\sum_{\substack{\iota_1,\ldots,\iota_r\in \mathrm f\{0,1\}^n  ,\\
\, j_1,\ldots,j_r\in  \{1\ldots n\} }}
p^{[k]}{}_{\iota_1,\ldots,\iota_r}^{j_1,\ldots,j_r}\,
(R/Q)^{\iota_1}_{j_1}\ldots (R/Q)^{\iota_r}_{j_r} \right)Q_k. \tag{F2L'}\label{F2L'}\]
Etc. The exact statements are left to the reader.

In this way, we obtain that for FQ conform operations we have $((n-2)2^{n-1}+1)^r2^{n-1}n$
(satisfying \texttt{(vL)}:  $((n-2)2^{n-1}+1)^rn$) real coefficients on the $r$th expansion level $(r\geq1)$ to be chosen freely;
and that for FQ conform orthogonalization procedures satisfying \texttt{(vL)}, \texttt{(vC')} we have
$2(((n-2)2^{n-1}+1)^r-1)(1-2^{-n})+1$
 real coefficients on the $r$th expansion level $(r\geq1)$ to be chosen freely.
\end{point}
\begin{remark} One can also encode a formal FQ conform-operation $\Psi$  (with $n$ variables) by a scalar FQ operation $\Upsilon$ and a vectorial FQ operation $\Phi$ (with $n-1$ variables).
Indeed, we can set
\[\Psi(A_1,A_2,\ldots,A_n)=(\tilde B_1A_1,\tilde B_2A_1,\ldots,\tilde B_nA_1)\]
where
\[\Upsilon(A_2 A_1^{-1},\ldots,A_nA_1^{-1})=\tilde B_1\qquad\text{and}\qquad \Phi(A_2 A_1^{-1},\ldots,A_nA_1^{-1})=(\tilde B_2,\ldots,\tilde B_n).\]
For formal FQ conform-orthogonalization procedures it is better suited to take
\[\Psi(A_1,A_2,\ldots,A_n)=(\tilde B_1A_1,\tilde B_2\tilde B_1A_1,\ldots,\tilde B_n\tilde B_1A_1).\]
\end{remark}

\section{The formal generalization of the symmetric  procedure}
is suggested by Corollary \ref{cor:charGS}: we should use the symmetric connection condition.
\begin{point} Let $\omega$ be a connection type as in \ref{po:connec}.
If $A=(A_1,\ldots,A_n)$, $Q=(Q_1,\ldots,Q_n)$, and $Q$ is a (floating) Clifford system, then
we define the associated iterators by
\[\Step_A^{\omega} Q:=\left(\Ad\exp\Pi^{\omega}_Q(A-Q) \right)Q;\]
and similarly in the floating case. They yield new (floating) Clifford systems.
\end{point}
\begin{remark} Throughout this paper, we use $\exp X$ is the usual sense.
But not much would change using $\exp^{\mathrm l} X=1+X$ or $\exp^{\mathrm r} X=(1-X)^{-1}$.
\end{remark}
\begin{theorem}\label{thm:iter1}
(a) Assume $A=(A_1,\ldots,A_n)$, $Q=(Q_1,\ldots,Q_n)$, and $Q$ is a (floating) Clifford system,
$A-Q\in(\mathfrak A^{(1)})^n$. Consider
\[Q^{[k]}=\left(\Step_A^{\omega}\right)^kQ. \]
Then, we claim,
\[Q^{[\omega]}=\lim_{k\rightarrow \infty}Q^{[k]}\]
exists; and in fact,
\[Q^{[k]}-Q^{[\omega]}\in (\mathfrak A^{(k+1)})^n;\]
furthermore, $Q^{[\omega]}$ is a Clifford system.
Moreover, it is the only Clifford system $\widetilde{Q}$ in $\mathfrak A$ such that
$\Pi^{\omega}_{\widetilde{Q}}A=0$
and $A-\widetilde{Q}\in(\mathfrak A^{(1)})^n$ hold.

(b) There is a corresponding statement for floating Clifford systems $Q$, with $\omega$
replaced by $\mathrm f\omega$, etc.
\begin{proof}
(a) First, we remark that due to the nature of conjugation, expansion orders are the same with respect for all $Q^{[k]}$.
From Lemma \ref{lem:formunic}.(i), it follows that
\[ \Pi^\omega_{\tilde Q}A=O( (R/Q)^{\geq r}\cdot Q)\quad \Rightarrow\quad\Pi^\omega_{(\Step^\omega_A\tilde Q)}A=O( (R/Q)^{\geq r+1}\cdot Q)\]
(the big ordo notation used reasonably). Using this observation, by induction, we can see that
\[\Pi^\omega_{ Q^{[k]}}A=O( (R/Q)^{\geq k+1}\cdot Q)\]
and
\[Q^{[k]}-Q^{[k+1]}\in (\mathfrak A^{(k+1)})^n.\]
This establishes the properties of the limit, which is obviously a Clifford system, and satisfies
the characterizing properties listed in the theorem.
Unicity is a consequence of Lemma \ref{lem:formunic}.(ii).
(b) The floating case is similar.
\end{proof}
\end{theorem}
\begin{lemma}\label{lem:formunic}

(a) Assume that $Q\in\Gr\mathfrak A^n$, $A\in \mathfrak A^n$ such that
$A- Q\in(\mathfrak A^{(1)})^n$, and $r\geq1$. Suppose that $\omega$ is a natural connection as in \ref{po:connec}.

(i)
If $X\in\mathfrak A^{(r)}\cap  T_{Q}\Gr\mathfrak A^n$ and $\widetilde{Q}=(\Ad\exp X)Q$, then
\[ X=\Pi_{Q}A-\Pi_{\widetilde{Q}}A\mod(\mathfrak A^{(r+1)})^n.\]

(ii) If $\widetilde{Q}\in\Gr\mathfrak A^n$, such that
$\widetilde{Q}- Q\in(\mathfrak A^{(r)})^n$, $(r\geq1)$, then
\[ \widetilde{Q}-Q=(\aad Q)\Pi_{Q}A-(\aad {\widetilde{Q}})\Pi_{\widetilde{Q}}A\mod(\mathfrak A^{(r+1)})^n.\]

(b) There are analogous statements in the floating case with $\Gr^{\mathrm f}$, $\Pi^{\mathrm f\omega}$, $\ad^{\mathrm f}$.
\begin{proof} We can pass from $\mathfrak A$ to $\mathfrak A/\mathfrak A^{(r+1)}$.
Then the basic observation is that $\mathfrak A^{(r)}$ and $\mathfrak A^{(1)}$ annihilate each other (multiplicatively). Then:

(a) We see that, in this setting, $Y\in\mathfrak A^{(1)}$ implies $(\Ad \exp X)Y=Y$;
$(\Ad \exp X)A=A+(\ad X)Q$; and $\widetilde{Q}=Q+(\ad X)Q$.
Then, due to the observations above, the connection identity \eqref{C}, and the conjugation-invariance \eqref{CN2},  we find
\begin{multline*}
\Pi_{\widetilde{Q}}A=\Pi_{(\Ad \exp X)Q}A= (\Ad \exp X)(\Pi_{Q}(\Ad \exp X)^{-1}A)= \\
=(\Ad \exp X)(\Pi_{Q}(A-(\ad X)Q))= (\Ad \exp X)(\Pi_{Q}A-X)=\Pi_{Q}A-X.
\end{multline*}
(b)
It can be assumed that $\widetilde{Q}$  is as in (i).
Then, in the same manner,
\begin{multline*}
(\aad \widetilde{Q})\Pi_{\widetilde{Q}}A=\ldots=(\Ad \exp X)(\aad Q)(\Pi_{Q}A-X)=(\Ad \exp X)(\aad Q)\Pi_{Q}A-(\aad Q)X)=\\
=(\Ad \exp X)(\aad Q)\Pi_{Q}A-(\widetilde{Q}-Q))=(\aad Q)\Pi_{Q}A-\widetilde{Q}+Q
\end{multline*}
shows the statement. (b) The proof is analogous in the floating case.
\end{proof}
\end{lemma}
\begin{defi}
We call the formal FQ orthogonalization procedure obtained above
\[\underline{\mathcal O}^{\omega}(A):=`` \left(\Step_A^{\omega}\right)^\infty Q  "\]
as the formal $\omega$-orthogonalization procedure;  similarly for
the formal FQ conform-ortho\-go\-na\-li\-za\-tion procedure with $\mathrm f\omega$.
\end{defi}
\begin{theorem}
(a) $\underline{\mathcal O}^{\omega}$ is a formal FQ orthogonalization procedure
(\texttt{(vR)},  \texttt{(Frm)}, \texttt{(\underline{Ana})}/ \texttt{(Nat)}, \texttt{(CP)}) with \texttt{(vL)}, \texttt{(vC)}, \texttt{(\underline{H}\,$\mathtt{{}^{0}}$)},
and satisfies the additional condition\texttt{(\underline{FSt})}.

In the case $\omega=\mathrm{GS}$ it is the the formal restriction of $\mathcal O^{\mathrm{GS}}$, and
in the case $\omega=\mathrm{Sy}$ it has the properties
\texttt{(\underline{PSy})}, \texttt{($\mathtt\Sigma$)},  \texttt{(O)}.

(b) $\underline{\mathcal O}^{\mathrm f\omega}$ is a formal FQ conform-orthogonalization procedure
(\texttt{(vR)},  \texttt{(Biv')}, \texttt{(Frm')}, \texttt{(\underline{Ana'})}/\texttt{(Nat)}, \texttt{(CP')}), with \texttt{(vL)}, \texttt{(vC')},
and satisfies the additional condition \texttt{(\underline{FSt})}.

In the case $\omega=\mathrm{GS}$ it is the the formal restriction of $\mathcal O^{\mathrm{fGS}}$, and
in the case $\omega=\mathrm{Sy}$ it has the properties
\texttt{(\underline{PSy'})}, \texttt{($\mathtt\Sigma$)},  \texttt{(O)}.
\end{theorem}
\begin{theorem}
$\underline{\mathcal O}^{\omega}$ and $\underline{\mathcal O}^{\mathrm f\omega}$  can be characterized as follows:

(a) $\underline{\mathcal O}^{\omega}(A)=Q$ if and only if (i) $Q$ is a Clifford system, (ii) $\Pi^{\omega}_QA=0$, (iii) $A-Q\in(\mathfrak A^{(1)})^n$.

(b) $\underline{\mathcal O}^{\mathrm f\omega}(A)=Q$ if and only if (i) $Q$ is a floating Clifford system, (
ii) $\Pi^{\mathrm f\omega}_QA=1^{\mathrm f}$, (iii) $A-Q\in(\mathfrak A^{(1)})^n$.
\begin{proof}[Proofs] The formal FQ (conform-)orthogonalization properties and the characterization
follows from the previous discussion, and the rest follows from the characterization.
\end{proof}
\end{theorem}
(For a connection as in \ref{po:connec}(iv) or (iv') we loose  \texttt{(vL)}.)
\begin{remark}
It is reasonable to call the condition
$\sum_{i=1}^n [A_i,Q_i]=0$  (cf. Lemma \ref{lem:conn}.(ii)) as the metric trace commutativity, \texttt{(MTC)}, of $A$ and $Q$;
and the condition
$\sum_{i=1}^n Q_i^{-1}A_i=\sum_{i=1}^n A_iQ_i^{-1}=n\cdot1$ (cf. Lemma \ref{lem:conn}.(ii')) as the metric trace inverse property,
\texttt{(MTI)}, with respect to $A$ and $Q^{-1}$.
\end{remark}
\begin{remark} Let $w=(w_1,\ldots,w_n)$ a weight as in \ref{po:connec}.
Then ordinary and floating cases are connected by
\[(1).\underline{\mathcal O}^{(w_1,\ldots,w_n)}(A_1,\ldots,A_n)=
\lim_{u\rightarrow\infty}\underline{\mathcal O}^{\mathrm f (u,w_1,\ldots,w_n)}(1,A_1,\ldots,A_n)\]
meaning concatenation of sequences on the left side.
Such a statement also can be formulated (then without even taking limits) to more general connection data.
\end{remark}

\begin{example} Using the notation of Example \ref{exa:gs},
we find that for $\Psi=\underline{\mathcal O}^{\mathrm{Sy}}$ we have coefficient matrices
$\mathbf P^{[1]}_0=[1]$, $ \mathbf P^{[2]}_0=[1]$, and
\[\mathbf P^{[1]}_1=\left[ \begin {array}{cccccccc} 0&0&1&\frac12&0&0&0&\frac12\end {array} \right],\]
\[\mathbf P^{[2]}_1=\left[ \begin {array}{cccccccc} 0&0&0&\frac12&0&1&0&\frac12\end {array}\right],\]
\[\mathbf P^{[1]}_2=\frac18
\left[ \begin {array}{cccccccc}
0&0&-4&-1&0&0&0&-1\\
0&0&-2&-2&0&0&0&-2\\
-4&-2&4&2&0&-2&0&2\\
-1&-2&2&1&-1&2&0&1\\
0&0&0&-1&0&0&0&-1\\
0&0&2&-2&0&0&-2&2\\
0&0&0&0&0&-2&0&0\\
-1&-2&2&1&-1&-2&0&1
\end {array} \right],\]
\[\mathbf P^{[2]}_2=\frac18
\left[ \begin {array}{cccccccc}
0&0&0&-1&0&0&0&-1\\
0&0&-2&0&0&0&0&0\\
0&-2&0&2&0&2&0&-2\\
-1&0&-2&1&-1&2&-2&1\\
0&0&0&-1&0&-4&0&-1\\
0&0&-2&2&-4&4&-2&2\\
0&0&0&-2&0&-2&0&-2\\
-1&0&2&1&-1&2&-2&1
\end {array} \right]
.\]
Notice that $\varpi=(15)(27)(36)(48)$ generates the symmetry $p^{[1]}_{i_1, \ldots } = p^{[2]}_{\varpi(i_1), \ldots } $.

More generally, for $\Psi=\underline{\mathcal O}^{\mathrm{GS}(t)}$ we have coefficient matrices
$\mathbf P^{[1]}_0=[1]$, $ \mathbf P^{[2]}_0=[1]$, and
\[\mathbf P^{[1]}_1=\left[ \begin {array}{cccccccc} 0&0&1& \frac1{1+t}&0&0&0
&{\frac {t}{1+t}}\end {array} \right]
,\]
\[\mathbf P^{[2]}_1=\left[ \begin {array}{cccccccc} 0&0&0& \frac1{1+t}&0&1&0
&{\frac {t}{1+t}}\end {array} \right]
,\]
\[\mathbf P^{[1]}_2=\frac12
 \left[ \begin {array}{cccccccc} 0&0&-1&- \frac1{(1+t)^2}&0&0
&0&-{\frac {t}{ \left( 1+t \right) ^{2}}}\\
0&0&- \frac1{1+t}&- \frac1{1+t}&0&0&0&-{\frac {t}{1+t}}\\
-1&- \frac1{1+t}&1& \frac1{1+t}&0&-{\frac {t}{1+t}}&0&{\frac {t}{1+t}}\\
- \frac1{(1+t)^2}&- \frac1{1+t}& \frac1{1+t}& \frac1{(1+t)^2}&-{\frac {t}{ \left( 1+t \right) ^{2}}}
 &{\frac {t}{1+t}}&0&{\frac {t}{ \left( 1+t \right) ^{2}}}\\
0&0&0&-{\frac {t}{ \left( 1+t \right) ^{2}}}&0&0&0&-{\frac {{t}^{2}}{ \left( 1+t \right) ^{2}}}\\
0&0&{\frac {t}{1+t}}&-{\frac {t}{1+t}}&0&0&-{\frac {t}{1+t}}&{\frac {t}{1+t}}\\
0&0&0&0&0&-{\frac {t}{1+t}}&0&0\\
-{\frac {t}{ \left( 1+t \right) ^{2}}}&-{\frac {t}{1+t}}&{\frac {t}{1+t}}&{\frac {t}{ \left(1+t \right) ^{2}}}
&-{\frac {{t}^{2}}{ \left( 1+t \right) ^{2}}}&-{\frac {t}{1+t}}&0&{\frac {{t}^{2}}{ \left( 1+t \right) ^{2}}}
\end {array} \right],
\]
\[\mathbf P^{[2]}_2=\frac12
\left[ \begin {array}{cccccccc}
0&0&0&- \frac1{(1+t)^2}&0&0&0&-{\frac {t}{ \left( 1+t \right) ^{2}}}\\
0&0&- \frac1{1+t}&0&0&0&0&0\\
 0&- \frac1{1+t}&0& \frac1{1+t}&0& \frac1{1+t}&0&- \frac1{1+t}\\
- \frac1{(1+t)^2}&0&- \frac1{1+t}& \frac1{(1+t)^2}&
-{\frac {t}{ \left( 1+t \right) ^{2}}}& \frac1{1+t}&- \frac1{1+t}&{\frac {t}{ \left( 1+t \right) ^{2}}}\\
0&0&0&-{\frac {t}{ \left( 1+t \right) ^{2}}}&0&-1
&0&-{\frac {{t}^{2}}{ \left( 1+t \right) ^{2}}}\\
0&0&- \frac1{1+t}& \frac1{1+t}&-1&1&-{\frac {t}
{1+t}}&{\frac {t}{1+t}}\\
0&0&0&- \frac1{1+t}&0&-{\frac {t}{1+t}}&0&-{\frac {t}{1+t}}\\
-{\frac {t}{ \left( 1+t \right) ^{2}}}&0& \frac1{1+t}&{\frac {t}{ \left( 1+t \right) ^{2}}}&-{\frac {{t}^{
2}}{ \left( 1+t \right) ^{2}}}&{\frac {t}{1+t}}&-{\frac {t}{1+t}}&{
\frac {{t}^{2}}{ \left( 1+t \right) ^{2}}}\end {array} \right]
.\]
This connects the Gram-Schmidt ($t=0$) and symmetric $(t=1)$ cases.
Notice that $\varpi=(15)(27)(36)(48)$ generates the symmetry $p^{[1]}_{i_1, \ldots }(t) = p^{[2]}_{\varpi(i_1), \ldots }(\frac1t)$.
\end{example}

\begin{theorem}
(a) Among the formal FQ orthogonalization  procedures $\Psi$ with \texttt{(vL)} and \texttt{(vC)},

(i) $\underline{\mathcal O}^{\mathrm{GS}}$ is  characterized by  \texttt{(Fil)} and \texttt{(\underline{FSt})};

(ii) $\underline{\mathcal O}^{\mathrm{Sy}}$ is  characterized by  \texttt{($\mathtt{\Sigma}$)} and \texttt{(\underline{FSt})}.
\\
(b) Among the formal FQ conform-orthogonalization  procedures $\Psi$ with \texttt{(vL)} and \texttt{(vC')},

(i') $\underline{\mathcal O}^{\mathrm{fGS}}$ is  characterized by  \texttt{(Fil)} and \texttt{(\underline{FSt})};

(ii') $\underline{\mathcal O}^{\mathrm{fSy}}$ is  characterized by  \texttt{($\mathtt{\Sigma}$)}  and \texttt{(\underline{FSt})}.

\begin{proof}
(i) It is not hard so prove that  \texttt{(Fil)} implies that the power series expansions of
$\Psi$ and $\underline{\mathcal O}^{\mathrm{GS}}$ are the the same up to first order.
In particular, if $Q=\Psi(A)$, $A=Q+R$, then $\Pi^{\mathrm{GS}}_Q(Q+R)=O( (R/Q)^{\geq2}\cdot Q)$.
Considering $\Psi(Q+tR) =Q$ (i.~e. \texttt{(\underline{FSt})}),
we obtain $\Pi^{\mathrm{GS}}_Q(Q+tR)=O( (tR/Q)^{\geq2}\cdot Q)$.
``Taking derivative in $t=0$'' we find $\Pi^{\mathrm{GS}}_QR=0$, which, by the previous theorem,  yields
$\underline{\mathcal O}^{\mathrm{GS}}(A)=Q$. (ii)-- The other cases are analogous.
\end{proof}
\end{theorem}

\section{The analytic extension of the symmetric  procedure}
We discuss two different methods here.
The first, more general one is the GS anchoring method;
the second, more natural one is the GS connection deformation method.
\begin{point}
Although we constructed $\underline{\mathcal O}^{\mathrm{Sy}}$ formally,
the nature of the recursion method in its construction implies,
by standard generating function arguments, that the power series expansion is analytic.
More precisely, if $Q\in\Gr\mathfrak A^n$, then
\[A\mapsto \text{``$\underline{\mathcal O}^{\mathrm{Sy}}_{\text{expanded around }Q}$''}(A)\]
extends $\underline{\mathcal O}^{\mathrm{Sy}}$ in a $0$-neighbourhood of $Q$, and eventually of $\Gr\mathfrak A^n$.
This is very nice, but, unfortunately, it fails short of producing a natural FQ
orthogonalization procedure, essentially because we have no control over $Q$ (and convergence).
This is corrected in the following
\end{point}
\begin{defi}
 Suppose that $\Psi$ is a formal FQ operation. Let
\[\partial_r\Psi(Q;R):= \frac{\partial^r}{\partial t^r}\Psi(Q+tR)\Bigr|_{t=0}\]
in terms of formal power series, where $t$ is a formal commutative variable.
(This is an $r$-homogenous polynomial in the variables $(R/Q)^\iota_j$.)

 If $P:[0,1]\rightarrow\mathfrak A^n$ is a real analytic function such that
\[ \frac{\partial^rP(t)}{\partial t^r}\Bigr|_{t=0}=\partial_r\Psi(\mathcal O^{\mathrm{GS}}(A) ;A-\mathcal O^{\mathrm{GS}}(A)),\]
then we set
\[\Psi_{\mathcal O^{\mathrm{GS}}}(A):=P(1).\]
This yields the Gram-Schmidt anchoring $\Psi_{\mathcal O^{\mathrm{GS}}}$ of $\Psi$.

We define the orthogonal Gram-Schmidt anchoring $\Psi_{\mathcal O^{\mathrm{GS}}/\mathrm O}(A)$ as the common value (if it exists) of
all well-defined quantities
\[ U^{-1}\Psi_{\mathcal O^{\mathrm{GS}}}(UA), \]
where $U\in\mathrm O(n)$ is an orthogonal change of coordinates.

Of course, we can define anchoring to any other FQ orthogonalization procedure, but at the moment we have only $\mathcal O^{\mathrm{GS}}$.
\end{defi}
\begin{theorem}
$\underline{\mathcal O}^{\mathrm{Sy}}_{\mathcal O^{\mathrm{GS}}/\mathrm O}$ is an FQ orthogonalization procedure
(\texttt{(vR)}, \texttt{(Nat)}, \texttt{(Ana)}, \texttt{(CP)}) with \texttt{(vL)}, \texttt{(vC)}, \texttt{(H\,$\mathtt{{}^{0}}$)}, which
generalizes the symmetric orthogonalization procedure (\texttt{(PSy)}, \texttt{($\mathtt\Sigma$)}, \texttt{(O)}),
and satisfies the additional condition \texttt{(FSt)}.
\begin{proof}
All the properties are immediate from the construction (via analytic continuation) except
the domain condition part of \texttt{(H\,$\mathtt{{}^{0}}$)}.
This however, follows from the following observation:
If the real analytic function $P:[0,1]\rightarrow\mathfrak A^n$ extends
\[t\mapsto \text{``$\underline{\mathcal O}^{\mathrm{Sy}}_{\text{expanded around }\mathcal O^{\mathrm{GS}}(A)}$''}
((1-t)\mathcal O^{\mathrm{GS}}(A)+tA)\]
for small $t$, then for $s>0$, the real analytic function \[P_s:[0,1]\rightarrow\mathfrak A^n\]
\[t\mapsto P\left(\frac{st}{(1-t)+st}\right)\]
extends
\[t\mapsto \text{``$\underline{\mathcal O}^{\mathrm{Sy}}_{\text{expanded around }\mathcal O^{\mathrm{GS}}(A)}$''}
((1-t)\mathcal O^{\mathrm{GS}}(A)+stA)\]
for small $t$ (as we know homogeneity for $(1-t)+st\approx1$ from  \texttt{(\underline{H}\,$\mathtt{{}^{0}}$)}).
\end{proof}
\end{theorem}
\begin{defi} Assume that $A=Q+R$, in the familiar formal setting.
 Let
\[\underline{\Omega}(A,t):= \underline{\mathcal O}^{\mathrm{GS}(t)}(Q+R)\]
understood such that $t$ is a formal variable commuting with everything.
This is means that terms like $(1+t^2+t^3)^{-1}$ are to be replaced by $1-(t^2+t^3)+(t^2+t^3)^2-\ldots$, etc.
Let
\[\underline{\Omega}_r(A):= \frac{\partial^r}{\partial t^r}\underline{\Omega}(A,t)\Bigr|_{t=0}.\]
\end{defi}
\begin{theorem}
Assume $A=(A_1,\ldots,A_n)$, $Q=(Q_1,\ldots,Q_n)=\mathcal O^{\mathrm{GS}}(A)$.
Let \[Q^{[0]}:=Q.\]
Then we define, recursively
\[X^{[k]}:=\partial\mathcal O^{\mathrm{GS}}( Q;\Pi_{Q^{[k]}}^{\mathrm{GS}(t)}(A-Q) ) ; \]
and
\[Q^{[k+1]}:= (\Ad \exp X^{[k]}) Q^{[k]} . \]
Then, we claim,
\[ \Pi_{Q^{[k]}}^{\mathrm{GS}(t)}(A-Q)  \in (t^{k+1} \mathfrak A)^n,\]
and
\[X^{[k]} \in t^{k+1} \mathfrak A.\]
Hence, eventually,
\[Q^{[k]}-  \underline{\Omega}(Q,t)  \in (t^{k+1} \mathfrak A)^n.\]
\begin{proof}
This is similar to the proof of Theorem \ref{thm:iter1}.
The crucial point is following: we can write $\Pi^{\mathrm{GS(t)}}_QA=\Pi^{\mathrm{GS}}_QA+t\Lambda^{t}_QA$, where
$\Lambda^{t}_QA$ is an appropriate formal expression.
Hence, when we make correction in $Q^{[k]}$ of order $t^{k+1}$ in order to make sure that
$\Pi^{\mathrm{GS(t)}}_{Q^{[k+1]}}A$ vanishes in the $t^{k+1}$ terms, it is sufficient to make the correction
as if making a variation in $\Pi^{\mathrm{GS}}_QA$, because the variation $t\Lambda^{t}_QA$ is of order $t^{k+2}$ anyway.
\end{proof}
\end{theorem}

\begin{theorem}
The expression $\text {``$\underline{\Omega}_r{}_{\text{expanded around  $\mathcal O^{\mathrm{GS}}(A)$}}$''}(A)$
can be extended analitically (and uniquely) as $\Omega_r(A)$ from a neighboorhood of $\Gr\mathfrak A^n$ to the domain of $\mathcal O^{\mathrm{GS}}$
along, say, straight segments $[\mathcal O^{\mathrm{GS}}(A),A]$ .

Moreover, $\Omega_r(A)$ can be realized as a finite integral formula using terms of $A$ and $Q=\mathcal O^{\mathrm{GS}}(A)$.
In particular, an extension of $\Omega_r(A)$ by a finite integral formula like above does not depend on the integral formula.
\begin{proof}
According the statement of the previous theorem $\underline{\Omega}_r(A)$ is realized as a finite expression
of $Q$, $A$, and $\partial\mathcal O^{\mathrm{GS}}( Q;-)$; hence a finite integral formula.
This, however provides an analytic extension on the segment, which is unique.
\end{proof}
\end{theorem}
\begin{defi}
 If $P:[0,1]\rightarrow\mathfrak A^n$ is a real analytic function such that
\[ \frac{\partial^r P(t)}{\partial t^r}\Bigr|_{t=0}=\Omega_r(A),\]
then we set
\[\mathcal O^{\mathrm{Sy}}_{\mathrm{GS}}(A):=P(1).\]
This yields the symmetric Gram-Schmidt connection deformation realization $\mathcal O^{\mathrm{Sy}}_{\mathrm{GS}}$.

We define the orthogonal Gram-Schmidt
connection deformation realization $\mathcal O^{\mathrm{Sy}}_{\mathrm{GS}/\mathrm O}$
as the common value (if it exists) of all well-defined quantities
\[ U^{-1}\mathcal O^{\mathrm{Sy}}_{\mathrm{GS}} (UA) \]
where $U\in\mathrm O(n)$ is an orthogonal change of coordinates.
\end{defi}
\begin{theorem}
$\underline{\mathcal O}^{\mathrm{Sy}}_{\mathrm{GS}/\mathrm O}$ is an FQ orthogonalization procedure
(\texttt{(vR)}, \texttt{(Nat)}, \texttt{(Ana)}, \texttt{(CP)}) with \texttt{(vL)}, \texttt{(vC)}, \texttt{(H\,$\mathtt{{}^{0}}$)}, which
generalizes the symmetric orthogonalization procedure (\texttt{(PSy)}, \texttt{($\mathtt\Sigma$)}, \texttt{(O)}),
and satisfies the additional condition \texttt{(FSt)}.
\begin{proof}
All the properties are immediate from the construction (via analytic continuation).
Here even the domain condition part of \texttt{(H\,$\mathtt{{}^{0}}$)} is easy, because
one can prove that  $\Omega_r$ itself is homogeneous.
\end{proof}
\end{theorem}
\begin{point}
One can also generalize the results of this section to the floating case.
The exact statements are left to the reader.

The GS anchoring method is easier to compute, while the GS connection deformation method
is more natural, because the domain condition seems to be more reasonable.

The are many other ways to extend the definition of the symmetric orthogonalization procedure,
some of them are with more natural invariance properties, but not necessarily better domain conditions than
in the case of the GS connection deformation method.
\end{point}

\begin{point}
The ideal solution would be having closed formulas as (hopefully canonical) realizations of FQ (conform-)or\-tho\-go\-na\-li\-za\-tion procedures.

In (very) small dimensions this is possible:

(i) The case of floating Clifford systems, $n=1$ is trivial: it yields the
identity operation $\mathcal O^{\mathrm{fSy}}((A_1))=(A_1)$.

(ii) The case of  Clifford systems, $n=1$ is well-known: it yields the
polarization operation $\mathcal O^{\mathrm{Sy}}((A_1))=(\pol A_1)$.

(iii) In the case of floating Clifford systems, $n=2$, we can take the definition
\[\mathcal O^{\mathrm{fSy}}(A_1,A_2):=\biggl(\frac{A_1+(\pol A_1A_2^{-1})A_2}{2},\frac{A_2+(\pol A_2A_1^{-1})A_1}{2} \biggr).\]

As the polarization operation satisfies the identities
$(\pol AB^{-1})B=-(\pol BA^{-1})B=-B(\pol A^{-1}B)=B(\pol B^{-1}A)$, we see that
\[\mathcal O^{\mathrm{fSy}}(A_1,A_2)=\frac12\biggl( (A_1,A_2)+ (-A_2,A_1)\cdot\pol A_1^{-1}A_2  \biggr).\]
Moreover, regarding the inverse floating system,
\[\lambda_1\mathcal O^{\mathrm{fSy}}(A_1,A_2)_1^{-1}+\lambda_2\mathcal O^{\mathrm{fSy}}(A_1,A_2)_2^{-1}=
\int_{t=0}^{2\pi}\frac{\lambda_1\cos t+\lambda_2\sin t}{A_1\cos t+A_2\sin t}\,\frac{\mathrm dt}{2\pi}\]
holds, in which form rotation invariance is even more transparent.

Notice that the domain condition is very simple and natural, the requirement is that
$(\lambda_1A_1+\lambda_2A_2)^{-1}$ should exist for $(\lambda_1,\lambda_2)\neq 0$.
This and some other nice properties make clear that this is \textit{the} right definition.

(iv) In the other cases, however, it is already not clear how to proceed.
\end{point}
This leads to

\section{Further directions}
There are, obviously, many questions left to pursue. For example:

(i) Regarding  the analytic extensions of FQ (conform-)orthonalization procedures constructed above or other ones:
One should obtain more precise information about their  domain;
find spectral characterizations of the procedures; or find explicit formulas for them.
More generally, one should study other, non-\texttt{(FSt)} FQ (conform-)orthonalization procedures, as
there are plenty of them.

(ii) Also, one should study general (and not necessarily vectorial) FQ operations:
One should get general statements regarding their invariance and arithmetic properties.
In particular, one could look for affine types, like a quantum inverse (dual) with sufficiently nice properties.

(iii) One should consider the limit $n=\infty$; in particular, using the continua $[0,1]$ and $\mathbb S^1$ for indexing;
hence allowing to take ``differentials in space-time''.

(iv) One should work out the bosonic case. Here the decompositions of the perturbations with respect
to the base systems will not be discrete anymore.\\~

In this way one hopes to get a fuller picture of noncommutative linear algebra (i--ii),
or, rather, functional analysis (iii--iv).

\section{Comparison to other approaches}
The study FQ operations can also be interpreted as
noncommutative analysis, or noncommutative functional calculus; where the slight difference in
emphasis is that in the first case one looks for an appropriate analytic environment to
deal with noncommutative objects (in our case: perturbations of Clifford systems),
while in the second case one wants to prolong commutative expressions (in our case: the classical orthogonalization procedures).
Work on these intertwined fields typically combines one or more of the following ideas:

(i) Extension of holomorphic function theory, in particular:
Utilizing noncommutative power series;
examining differentials, using nilpotent or Grassmannian variables;
considering compatibility with matrix calculus;
defining functional calculus around (generalized) spectra (cf.~%\cite{Ta0},
\cite{Tay}, \cite{Voi}, \cite{KVV});
extending the Cauchy formula (cf.~\cite{BDS},\cite{DSS}).

(ii) Drawing conclusions, axiomatizing, or just simply examining properties of concrete cases, in particular:
Doing any nonlinear computation with noncommuting algebra elements; or,
reinterpreting nonlinear expressions of noncommuting actions on  model spaces.

(iii) Making physical computations from quantum theory rigorous (cf.~\cite{Sim}, \cite{GJ}, \cite{BSZ}).
In particular, expounding on special features like:

(iv) Considering compact or, rather, nuclear perturbations of basic systems.

(v) Taking time time-ordered products
in discrete commuting patches (cf.~\cite{Mas}, \cite{NSS}, \cite{Jef});
or, along continua (often related to the Wiener measure), (cf.~\cite{Fey}, \cite{JL}).

(vi) Using generalized measure spaces and special transforms (cf.~\cite{BSZ}, \cite{AHM}, \cite{HKPS}, \cite{Kap}).

(vii) Building on bosonic quantites, or the Weyl algebra (cf.~\cite{Wey}, \cite{Fol}).

(viii) Building on fermionic quantities, or the Clifford algebra
(cf.~\cite{Jef}, \cite{CSS}, \cite{Alp}).

Often, one has partial results, but with an eye toward the bigger picture.
A general pattern is, however, that when we are strong on the constructive side,
we remain weak in terms of nonlinear arithmetic (like typical functional calculi);
and conversely,  in cases when good arithmetic properties, computational relevancy, are expected (like physics), constructions are weak.
In particular, it is not clear what a successful noncommutative functional calculus should look like even in the case of finite matrices.

Indeed, when we  put our present setup in the context above, it turns out to be quite simple but nontrivial:
We have worked with noncommutative power series enriched with Clifford variables (fermionic case)
which makes the structure richer, in fact, allows nontrivial analysis, but avoids analytic complications (opposed to the bosonic case).
Taken in a narrow sense,  we deal with finite-dimensional matrices.
We can see that it is not always easy to transfer ideas from the commutative case to the noncommutative one, but sometimes  possible.
So, one hopes that this setting might help to gather more information regarding arithmetically strong functional calculi
on finite matrices and hence in other cases.

%\section{Acknowledgements}

\end{document}